\documentclass[preprint,12pt]{elsarticle}
\usepackage{amsfonts}
\usepackage{amsmath}
\usepackage{amssymb}

\input{tcilatex}

\begin{document}

\begin{frontmatter}



\title{
New results on systems of generalized vector quasi-equilibrium problems 
}

\author{Monica Patriche}

\address{
University of Bucharest
Faculty of Mathematics and Computer Science
    
14 Academiei Street
   
 010014 Bucharest, 
Romania
    
monica.patriche@yahoo.com }

\begin{abstract}
In this paper, we firstly prove the existence of the equilibrium for the gene- ralized abstract economy. 
We apply these results to show the existence of solutions for systems of vector quasi-equilibrium 
problems with multiva- lued trifunctions. Secondly, we consider the generalized strong vector 
quasi-equilibrium problems and study the existence of their solutions in the case when the correspon
dences are weakly naturally quasi-concave or weakly biconvex and also in the case of weak-continuity 
assumptions. In all situations, fixed-point theorems are used.
\end{abstract}


\end{frontmatter}

% \linenumbers


\label{}



%% References


\bibliographystyle{elsarticle-num}
\bibliography{<your-bib-database>}

\begin{thebibliography}{99}
\bibitem{ans2} Q.H. Ansari, J.C. Yao, An existence result for the
generalized vector equilibrium, Appl. Math. Lett. 12 (1999) 53--56.

\bibitem{aum} R. Aumann, S. Hart, \textit{Bi-convexity and bi-martingales,}
Isr J Math \textbf{54, }\textit{2} (1986), 159-180.

\bibitem{bor} A. Borglin and H. Keiding, \textit{Existence of equilibrium
actions and of equilibrium: A note on the 'new' existence theorem,} J. Math.
Econom. \textbf{3} (1976), 313-316.

\bibitem{cai} G.Cai, S. Bu, Strong and weak convergence theorems for general
mixed equilibrium problems and variational inequality problems and fixed
point problems in Hilbert spaces, Journal of Computational and Applied
Mathematics, 247, \textit{1} (2013), 34-52.

\bibitem{deb} G. Debreu, \textit{A social equilibrium existence theorem.}
Proc. Nat. Acad. Sci. \textbf{38} (1952), 886-893.

\bibitem{ding} X. P. Ding and G. X. Z. Yuan, \textit{The study of existence
of equilibria for generalized games without lower semicontinuity} \textit{in
localy topological vector spaces,} J. Math. Anal. Appl. 227 (1998), 420-438.

\bibitem{ding2} X. Ding and He Yiran, \textit{Best Approximation Theorem for
Set-valued Mappings without Convex Values and} \textit{Continuity,} Appl
Math. and Mech. English Edition, 19, \textit{9} (1998), 831-836.

\bibitem{fan} K. Fan, \textit{Fixed-point and minimax theorems in locally
convex topological linear spaces,} Proc. Nat. Acad. Sci. U.S.A. 38(1952),
121-126.

\bibitem{gong} X.H. Gong, Strong vector equilibrium problem, J. Global
Optim. 36 (2006) 339--349.

\bibitem{gor} J. Gorski, F. Pfeuffer, K. Klamroth, \textit{Biconvex sets and
optimization with biconvex functions: a survey and extension, }Math. Meth.
Oper. Res. 66 (2007), 373-407.

\bibitem{he} H. He, S. Liu, Y. J. Cho, An explicit method for systems of
equilibrium problems and fixed points of infinite family of nonexpansive
mappings, Journal of Computational and Applied Mathematics, 235, \textit{14}%
, (2011), 4128-4139.

\bibitem{him} C. J. Himmelberg (1972), Fixed points of compact
multifunctions. J. Math. Anal. Appl. 38:205-207.

\bibitem{hou} S.H. Hou, X.H. Gong, X.M. Yang, Existence and stability of
solutions for generalized strong vector equilibrium problems with
trifunctions, J. Optim. Theory Appl. 146, \textit{2}, (2010) 387-398.

\bibitem{kim} W. K. Kim and K. K. Tan, New existence theorems of equilibrium
and applications Nonlinear Anal., 47 (2001), 531-542

\bibitem{lin} Y. C. Lin, On generalized vector equilibrium problems,
Nonlinear Analysis 70 (2009) 1040--1048.

\bibitem{lin2} L. Lin, Z. Yu, Q. Ansari and L. Lai, Fixed point and maximal
element theorems with applications to abstract economies and minimax
inequalities, J. of Math. Anal. Appl. 284, \textit{2}, (2003), 656--671.

\bibitem{lin3} L.J. Lin, L.F. Chen and Q.H. Ansari, Generalized abstract
economy and systems of generalized vector quasi-equilibrium problems, J.
Comput. Appl. Math 208 (2007), 341-352.

\bibitem{liu} Q.-M. Liu, L. Fan, G. Wang, Generalized vector
quasi-equilibrium problems with set-valued mappings. Appl. Math. Lett., 21
(9), (2008), 946-950

\bibitem{long} X.J. Long, N. J. Huang and K. L. Teo, Existence and stability
of solutions for generalized strong vector quasi-equilibrium problem, Math.
Comp. Model. 47 (2008), 445-451.

\bibitem{luc} D.T. Luc, Theory of Vector Optimization, in: Lecture Notes in
Economics and Mathematics Systems, vol. 319, Springer-Verlag, New York,1989.

\bibitem{orive} R. Orive, Z. Garc\'{\i}a, On a class of equilibrium problems
in the real axis, Journal of Computational and Applied Mathematics, 235, 
\textit{4}, (2010), 1065-1076.

\bibitem{pat} M. Patriche, Fixed point theorems and applications in theory
of games, Fixed point theory, to appear.

\bibitem{sha} Shafer W., Sonnenschein H.: Equilibrium in abstract economies
without ordered preferences. J. Math. Econom. \textbf{2,} 345-348 (1975).

\bibitem{tul} C. I. Tulcea, \textit{On the Approximation of Upper
Semi-continuous Correspondences and the Equilibriums of Generalized Games,}
J. Math. An. Appl. \textbf{136} (1988), 267-289.

\bibitem{yuan} X. Z. Yuan, \textit{The Study of Minimax Inequalities and
Applications to Economies and} \textit{Variational Inequalities, }Memoirs of
the American Society, 625(1988).

\bibitem{wang} Z. Wang, Y. Su, D. Wang, Y. Dong, A modified Halpern-type
iteration algorithm for a family of hemi-relatively nonexpansive mappings
and systems of equilibrium problems in Banach spaces, Journal of
Computational and Applied Mathematics, 235, \textit{8}, (2011), 2364-2371.
\end{thebibliography}






% \end{thebibliography}

\section{Introduction}

The vector equilibrium problem is a unified model of several problems, for
instance, vector variational inequalities, vector optimization problems or
Debreu-type equilibrium problems. For new results on this topic, the reader
is refered to \cite{ans2},\cite{cai},\cite{gong},\cite{he},\cite{hou},\cite%
{lin},\cite{lin3},\cite{liu}, \cite{long},\cite{orive}, \cite{wang}.

There are a lot of approaches in order to establish existence results for
solutions of vector equilibrium problems. Lin, Chen and Ansari \cite{lin3}
used the existence theorems for generalized abstract economies. This method
was proven to be a fruitful one, and it deserves to be deepened. The
generalized abstract economy was introduced by Kim and Tan \cite{kim} and
generalizes the previous models of abstract economies due to Debreu \cite%
{deb}, Shafer and Sonnenshine \cite{sha} or Borglin and Keiding \cite{bor}.
Kim and Tan motivated their work by the fact that any preference of a real
agent could be unstable by the fuzziness of consumers' behaviour or market
situations.

In this paper, by following the direction opened by Lin, Chen and Ansari 
\cite{lin3}, we firstly prove the existence of the equilibrium for the
generalized abstract economy with upper semicontinuous fuzzy
correspondences. We apply these results to show the existence of solutions
for systems of vector quasi-equilibrium problems with multivalued
trifunctions. Secondly, we consider the generalized strong vector
quasi-equilibrium problems and study the existence of their solutions in the
case when the correspondences are weakly naturally quasi-concave or weakly
biconvex and also in the case of weak-continuity assumptions. In both last
situations, fixed-point theorems are used.

The paper is organised as follows. Section 2 contains preliminaries and
notations. In Sections 3, the equilibrium existence of the generalized
abstract economy model is obtained. Section 4 studies the existence of
solutions for systmes of generalized vector quasi-equilibrium problems.
Section 5 presents types of convexity conditions which are sufficient in
order to guarantee that the generalized vector quasi-equilibrium problems
can be solved. The case of weak-continuity assumptions is approached at the
end.

\section{Preliminaries and notation}

For the reader's convenience, we present several properties of the
correspondences which are used in our proofs.

Let $X$ be a subset of a topological vector space $E$. The set $X$ is said t%
\textit{o have the property }$(K)$ if, for every compact subset $B$ of $X$,
the convex hull co$B$ is relatively compact in $E$. It is clear that each
compact convex set in a Hausdorff (resp., locally) topological vector space
always has property $(K)$. A normal topological space in wich each open set
is an $F_{\sigma }$ is called \textit{perfectly normal.}

Let $X$, $Y$ be topological spaces and $T:X\rightarrow 2^{Y}$ be a
correspondence. $T$ is said to be \textit{upper semicontinuous} if for each $%
x\in X$ and each open set $V$ in $Y$ with $T(x)\subset V$, there exists an
open neighborhood $U$ of $x$ in $X$ such that $T(x)\subset V$ for each $y\in
U$. $T$ is said to be \textit{lower semicontinuous} if for each x$\in X$ and
each open set $V$ in $Y$ with $T(x)\cap V\neq \emptyset $, there exists an
open neighborhood $U$ of $x$ in $X$ such that $T(y)\cap V\neq \emptyset $
for each $y\in U$. $T$ is said to have \textit{open lower sections} if $%
T^{-1}(y):=\{x\in X:y\in T(x)\}$ is open in $X$ for each $y\in Y.$ $T$ is
said to be \textit{compact} if, for any $x\in X,$ there exists an open
neighborhod $V(x)$ such that $T(N(x))=\cup _{y\in N(x)}T(y)$ is relatively
compact in Y.

The set valued map $\overline{T}$ is defined by $\overline{T}(x):=\{y\in
Y:(x,y)\in $cl$_{X\times Y}$ Gr $T\}$ (the set cl$_{X\times Y}$ Gr $(T)$ is
called the adherence of the graph of $T$)$.$ It is easy to see that cl $%
T(x)\subset \overline{T}(x)$ for each $x\in X.$ $T$ is said to be \textit{%
quasi-regular} if it has non-empty convex values, open lower sections and $%
\overline{T}(x)=$cl$T(x)$ for each $x\in X.$ $T$ is said to be \textit{%
regular} if it is quasi-regular and has an open graph.

If $X$ and $Y$ are topological vector spaces, $K$ is a non-empty subset of $%
X,$ $C$ is a non-empty closed convex cone and $T:K\rightarrow 2^{Y}$ is a
correspondence, then \cite{luc}, $T$ is called \textit{upper }$C$\textit{%
-continuous at} $x_{0}\in K$ if, for any neighbourhood $U$ of the origin in $%
Y,$ there is a neighbourhood $V$ of $x_{0}$ such that, for all $x\in V,$ $%
T(x)\subset T(x_{0})+U+C.$ $T$ is called \textit{lower }$C$\textit{%
-continuous at} $x_{0}\in K$ if, for any neighbourhood $U$ of the origin in $%
Y,$ there is a neighbourhood $V$ of $x_{0}$ such that, for all $x\in V,$ $%
T(x_{0})\subset T(x)+U-C.$

Now, we are presenting the approximation of upper semicontinuous
correspondences due to C. I. Tulcea

Let $X$ be a nonempty set, let $Y$ be a nonempty subset of a topological
vector space $E$, and let $T:X\rightarrow 2^{Y}.$ A family $(f_{j})_{j\in J}$
of correspondences between $X$ and $Y$, indexed by a nonempty filtering set $%
J$ (denote by $\leq $ the order relation in $J$), is an \textit{upper
approximating family for} $F$ \cite{tul} if (1) $T(x)\subset f_{j}(x)$ for
all $x\in X$ and all $j\in J;$ (2) for each $j\in J$ there is a $j^{\ast
}\in J$ such that, for each $h\geq j^{\ast }$ and $h\in J,$ $f_{h}(x)\subset
f_{j}(x)$ for each $x\in X$ and (3) for each $x\in X$ and V$\in $\ss , where 
\ss\ is a base for the zero neighborhood in $E$, there is a $j_{x,V}\in J$
such that $f_{h}(x)\subset T(x)+V$ if $h\in J$ and $j_{x,V}\leq h.$ From
(1)-(3), it is easy to deduce that that for each $x\in X$ and $k\in J,$ $%
T(x)\subset \cap _{j\in J}f_{j}(x)=\cap _{k\leq j,k\in J}f_{j}(x)\subset $cl$%
T(x)\subset \overline{T}(x).$

By observing Theorem 3 and the Remark of Tulcea [\cite{tul}, p.280 and pp
281-282], we have the following:

\begin{lemma}
(see\textbf{\ }\cite{ding}\textbf{).} \textit{Let }$(X_{i})_{i\in I}$\textit{%
\ be a family of paracompact spaces and let }$(Y_{i})_{i\in I}$\textit{\ be
a family of nonempty closed convex subsets, each in a locally convex
Hausdorff topological vector space and each has property }$(K)$\textit{. For
each }$i\in I,$ let\textit{\ }$T_{i}:X_{i}\rightarrow 2^{Y_{i}}$\textit{\ be
compact and upper semicontinuous with nonempty and convex values. Then,
there is a common filtering set }$J$\textit{\ (independent of }$i\in I$%
\textit{) such yhat, for each }$i\in I$\textit{, there is a family }$%
(f_{ij})_{j\in J}$\textit{\ of correspondences between }$X_{i}$\textit{\ and 
}$Y_{i\text{ }}$\textit{with the following properties:}
\end{lemma}

(i)\textit{\ for each }$j\in J,$ $(f_{ij})_{j\in J}$\textit{\ is regular;}

(ii)\textit{\ }$(f_{ij})_{j\in J}$\textit{\ and }$(\overline{f}_{ij})_{j\in
J}$ \textit{are upper approximating families for }$F_{i};$

(iii\textit{) for each }$j\in J,$ the correspondence $\overline{f}_{ij}$%
\textit{\ is continuous if }$Y_{i\text{ }}$\textit{\ is compact.}

Lemma 1 is a version of Lemma 1.1 in \cite{yuan} ( for $D=Y,$ we obtain
Lemma 1.1 in \cite{yuan}).

\begin{lemma}
(see \cite{pat}) Let $X$ be a topological space, $Y$ be a non-empty subset
of a locally convex topological vector space $E$ and $T:X\rightarrow 2^{Y}$
be a correspondence$.$ Let \ss\ be a basis of neighbourhoods of $0$ in $E$
consisting of open absolutely convex symmetric sets. Let $D$ be a compact
subset of $Y$. If for each $V\in $\ss , the correspondence $%
T^{V}:X\rightarrow 2^{Y}$ is defined by $T^{V}(x)=(T(x)+V)\cap D$ for each $%
x\in X,$ then $\cap _{V\in \text{\ss }}\overline{T^{V}}(x)\subseteq 
\overline{T}(x)$ for every $x\in X.$
\end{lemma}

The following lemma is an important result concerning the continuity of
correspondences which will be used in our proofs.

\begin{lemma}
(see \cite{yuan}). \textit{Let }$X$\textit{\ and }$Y$\textit{\ be two
topological spaces and let }$D$\textit{\ be an open subset of }$X.$\textit{\
Suppose }$T_{1}:X\rightarrow 2^{Y}$\textit{\ , }$T_{2}:X\rightarrow 2^{Y}$%
\textit{\ are upper semicontinuous correspondences such that }$%
T_{2}(x)\subset T_{1}(x)$\textit{\ for all }$x\in D.$\textit{\ Then the
correspondence }$T:X\rightarrow 2^{Y}$\textit{\ defined by}
\end{lemma}

\begin{center}
$T\mathit{(z)=}\left\{ 
\begin{array}{c}
T_{1}(x)\text{, \ \ \ \ \ \ \ if }x\notin D\text{, } \\ 
T_{2}(x)\text{, \ \ \ \ \ \ \ \ \ \ if }x\in D%
\end{array}%
\right. $
\end{center}

\textit{is also upper semicontinuous.\medskip }

The property of properly $C-$quasiconvexity for correspondences is presented
below.

Let $X$ be a non-empty convex subset of a topological vector space\textit{\ }%
$E,$ $Z$ be a real topological vector space, $Y$ be a subset of $Z$ and $C$
be a pointed closed convex cone in $Z$ with its interior int$C\neq \emptyset
.$ Let $T:X\rightarrow 2^{Z}$ be a correspondence with non-empty values. $T$
is said to be \textit{properly }$C-$\textit{quasiconvex on} $X$, iff for any 
$x_{1},x_{2}\in X$ and $\lambda \in \lbrack 0,1],$ either $T(x_{1})\subset
T(\lambda x_{1}+(1-\lambda )x_{2})+C$ or $T(x_{2})\subset T(\lambda
x_{1}+(1-\lambda )x_{2})+C.$

\section{\textbf{Equilibrium existence for generalized abstract economies}}

Because of the fuzziness of consumers' behaviour or market situations, in a
real market, any preference of a real agent would be unstable. Therefore,
Kim and Tan \cite{kim} introduced the fuzzy constraint correspondences in
defining the following generalized abstract economy.

Let $I$ be any set of agents (countable or uncountable). For each $i\in I$,
let $X_{i}$ be a nonempty set of actions available to the agent $i$ in a
topological vector space $E_{i}$ and $X=\underset{i\in I}{\prod }X_{i}.$

\begin{definition}
\cite{kim}A \emph{generalized abstract economy} $\Gamma
=(X_{i},A_{i},F_{i},P_{i})_{i\in I}$ is defined as a family of ordered
quadruples $(X_{i},A_{i},F_{i},P_{i})$ where $A_{i}:X\rightarrow 2^{X_{i}}$
is a constraint correspondence such that $A_{i}(x)$ is the state attainable
for the agent $i$ at $x$, $F_{i}:X\rightarrow 2^{X_{i}}$ is a fuzzy
constraint correspondence such that $F_{i}(x)$ is the unstable state for the
agent $i$ and $P_{i}:X\times X\rightarrow 2^{X_{i}}$ is a preference
correspondence such that $P_{i}(x,x)$ is the state preferred by the agent $i$
at $x.$
\end{definition}

\begin{definition}
An \emph{equilibrium }for $\Gamma $ is a point $(x^{\ast },y^{\ast })\in
X\times X$ such that for each $i\in I,$ $x_{i}^{\ast }\in $ $A_{i}(x^{\ast
}),$ $y_{i}^{\ast }\in $ $F_{i}(x^{\ast })$ and $P_{i}(x^{\ast },y^{\ast
})\cap A_{i}(x^{\ast })=\emptyset .$
\end{definition}

If for each $i\in I$ and each $x\in X,$ $F_{i}(x)=X_{i}$ and the preference
correspondence $P_{i}$ satisfies $P_{i}(x,y)=P_{i}(x,y^{^{\prime }})$ for
each $x,y,y^{^{\prime }}\in X,$ the definition of a generalized abstract
economy and an equilibrium coincide with the usual definitions of an
abstract economy and an equilibrium due to Shafer and Sonnenschein \cite{sha}%
.

The following theorem is the compact version of Theorem 5.1 of Lin and al. (%
\cite{lin2}). The set X is compact and the correspondences $A_{i}$, $F_{i}$
and $P_{i}$ have open lower sections.

\begin{theorem}
\textit{For each }$i\in I$\textit{\ (}$I$\textit{\ finite)}$,$\textit{\ let }%
$X_{i}$\textit{\ be a nonempty compact convex subset of a topological vector
space }$E_{i}$\textit{, }$X=\tprod\limits_{i\in I}X_{i},$\textit{\ }$%
A_{i}:X\rightarrow 2^{X_{i}}$\textit{\ a constraint correspondence, }$%
P_{i}:X\times X\rightarrow 2^{X_{i}}$\textit{\ a preference correspondence
and }$F_{i}:X\rightarrow 2^{X_{i}}$\textit{\ a fuzzy constraint
correspondence. Assume that the following conditions hold:}
\end{theorem}

(i)\textit{\ For all }$x\in X$\textit{, }$A_{i}(x)$\textit{\ and }$F_{i}(x)$%
\textit{\ are nonempty and convex;}

(ii)\textit{\ For all }$y_{i}\in X_{i},$\textit{\ }$A_{i}^{-1}(y_{i}),$%
\textit{\ }$F_{i}^{-1}(y_{i})$\textit{\ and }$P_{i}^{-1}(y_{i})$\textit{\
are open sets}$;$

(iii)\textit{\ For all (}$x,y)\in X\times X$\textit{, }$x_{i}\notin $co$%
P_{i}(x,y);$

(iv)\textit{\ The set }$W_{i}:$\textit{\ }$=\left\{ (x,y)\in X\times X:x_{i}%
\text{ }\in A_{i}(x)\text{ }and\text{ }y_{i}\text{ }\in F_{i}(x)\right\} $%
\textit{\ is closed in }$X\times X$\textit{;}

\textit{Then there exists }$(x^{\ast },y^{\ast })\in X$\textit{\ }$\times X$%
\textit{\ such that for each }$i\in I$\textit{, }$x_{i}^{\ast }\in
A_{i}(x^{\ast }),$\textit{\ }$y_{i}^{\ast }\in F_{i}(x^{\ast })$\textit{\
and }$A_{i}(x^{\ast })\cap P_{i}(x^{\ast },y^{\ast })=\emptyset $\textit{.
\medskip }

We state the following result which is an equilibrium existence theorem for
a generalized abstract economy with upper semicontinuous correspondences. We
use a method of approximation of upper semi-continuous correspondences
developed by C. I. Tulcea in \cite{tul}.

\begin{theorem}
\textit{For each }$i\in I$\textit{\ (}$I$\textit{\ finite)}$,$\textit{\ let }%
$X_{i}$\textit{\ be a nonempty compact convex subset with property (K) of a
topological vector space }$E_{i}$\textit{, }$X=\tprod\limits_{i\in I}X_{i}$%
\textit{\ be perfectly normal}$,$\textit{\ }$A_{i}:X\rightarrow 2^{X_{i}}$%
\textit{\ a constraint correspondence, }$P_{i}:X\times X\rightarrow
2^{X_{i}} $\textit{\ a preference correspondence and }$F_{i}:X\rightarrow
2^{X_{i}}$\textit{\ a fuzzy constraint correspondence. Assume that the
following conditions hold:}
\end{theorem}

(i)\textit{\ For all }$x\in X$\textit{, }$A_{i}(x)$\textit{\ and }$F_{i}(x)$%
\textit{\ are nonempty and convex;}

(ii)\textit{\ For all }$y_{i}\in X_{i},$\textit{\ }$P_{i}^{-1}(y_{i})$%
\textit{\ are open sets and the correspondences }$F_{i},A_{i}$\textit{\ are
upper semi-continuous, compact, with nonempty convex closed values}$;$

(iii)\textit{\ For all (}$x,y)\in X\times X$\textit{, }$x_{i}\notin $co$%
P_{i}(x,y);$

(iv)\textit{\ The set }$U_{i\text{ }}:=\{(x,y)\in X\times X:P_{i}(x,y)\cap
A_{i}(x)\neq \emptyset \}$\textit{\ is open.}

\textit{Then there exists }$(x^{\ast },y^{\ast })\in X$\textit{\ }$\times X$%
\textit{\ such that for each }$i\in I$\textit{, }$x_{i}^{\ast }\in \overline{%
A_{i}}(x^{\ast })$, \textit{\ }$y_{i}^{\ast }\in \overline{F_{i}}(x^{\ast })$%
\textit{\ and }$A_{i}(x^{\ast })\cap P_{i}(x^{\ast },y^{\ast })=\emptyset $%
\textit{. }

\textit{Proof.} According to Lemma 1, there is a common filtering set J such
that, for every $i\in I$, there exists a family $(A_{ij})_{j\in J}$ of
regular correspondences between $X$ and $X_{i},$ such that both $%
(A_{ij})_{j\in J}$ and $(\overline{A_{ij}})_{j\in J}$ are upper
approximating families for $A_{i}$ and a family $(F_{ij})_{j\in J}$ of
regular correspondences between $X$ and $X_{i},$ such that both $%
(F_{ij})_{j\in J}$ and $(\overline{F_{ij}})_{j\in J}$ are upper
approximating families for $F_{i}.$\ The correspondences $A_{ij}$ and $%
F_{ij} $ are regular, it follows that $A_{ij}$ and $F_{ij}$ have an open
graph and thus they have open lower sections. Since $A_{i}$ and $F_{i}$ have
closed graphs, the set $W_{i}:=\left\{ (x,y)\in X\times X:x_{i}\text{ }\in
A_{i}(x)\text{ and }y_{i}\text{ }\in F_{i}(x)\right\} $\textit{\ }is closed
in\textit{\ }$X\times X$\textit{. }Therefore the abstract economy $\Gamma
_{j}=(X_{i},A_{ij},P_{i},F_{ij})_{i\in I}$ satisfies all hypotheses of
Theorem 1. It follows by Theorem 1 that $\Gamma _{j}$ has an equilibrium $%
(x^{\ast j},y^{\ast j})\in X\times X$ such that $A_{ij}(x^{\ast j})\cap
P_{i}(x^{\ast j},y^{\ast j})=\emptyset ,$ $x_{i}^{\ast j}\in A_{ij}(x^{\ast
j})$ and $y_{i}^{\ast j}\in F_{ij}(x^{\ast j})$ for all $i\in I.$ Since $%
A_{i}(x^{\ast j})\subset A_{ij}(x^{\ast j}),$ it follows that $A_{i}(x^{\ast
j})\cap P_{i}(x^{\ast j},y^{\ast j})=\emptyset .$ Therefore $\{(x^{\ast
j},y^{\ast j})\}_{j\in I}\subset U_{i\text{ }}^{C}$ which is closed in $%
X\times X$ by condition (iv). On the other hand, note that $(x^{\ast
j},y^{\ast j})_{j\in I}$ is a net in the compact set $X\times X$; without
loss of generality, we may assume that $(x^{\ast j})_{j\in I}$ converges to $%
x^{\ast }\in X$ and $(y^{\ast j})_{j\in I}$ converges to $y^{\ast }\in X.$
Then, for each $i\in I,$ $x_{i}^{\ast }=\lim_{j\in I}x_{i}^{\ast j}$ and $%
y_{i}^{\ast }=\lim_{j\in I}y_{i}^{\ast j}.$ As $(x^{\ast },y^{\ast })\in U_{i%
\text{ }}^{C}$ for all $i\in I,$ it follows that $A_{i}(x^{\ast })\cap
P_{i}(x^{\ast },y^{\ast })=\emptyset $\textit{. }Since\textit{\ }$(x^{\ast
j},y^{\ast j})$ is an equilibrium point of $\Gamma _{j},$ we have that $%
x_{i}^{\ast j}\in A_{ij}(x^{\ast j})\subset \overline{A_{ij}}(x^{\ast j})$
and $y_{i}^{\ast j}\in F_{ij}(x^{\ast j})\subset \overline{F_{ij}}(x^{\ast
j}).$ As $\overline{A_{ij}}$ and $\overline{F_{ij}}$ have a closed graph, it
follows that $(x^{\ast },x_{i}^{\ast })\in $Gr$\overline{A_{ij}}$ and $%
(x^{\ast },y_{i}^{\ast })\in $Gr$\overline{A_{ij}}$ for every $i\in I.$ For
each $i\in I,$ since $(A_{ij})_{j\in I}$ is an upper approximation family
for $A_{i},$ it follows that $\cap _{j\in I}\overline{A_{ij}}(x)\subset 
\overline{A_{i}}(x)$ for each $\in X,$ so that $(x^{\ast },x_{i}^{\ast })\in 
$Gr$\overline{A_{i}}.$ Also, we have that$(x^{\ast },y_{i}^{\ast })\in $Gr$%
\overline{F_{i}}.$ Therefore, for each $i\in I,$ $A_{i}(x^{\ast })\cap
P_{i}(x^{\ast },y^{\ast })=\emptyset $, $x_{i}^{\ast }\in \overline{A_{i}}%
(x^{\ast })$ and\textit{\ }$y_{i}^{\ast }\in \overline{F_{i}}(x^{\ast }).$

The next results from this section prove the existence of equilibrium for
generalized abstract economies.

\begin{theorem}
\textit{Let }$\Gamma =\left\{ X_{i},A_{i},F_{i},P_{i}\right\} _{i\in I}$%
\textit{\ be a generalized abstract economy, where }$I$\textit{\ is any
index set, such that, for each }$i\in I$\textit{:}
\end{theorem}

i)\textit{\ }$X_{i}$\textit{\ is a non-empty convex subset of a Hausdorff
locally convex space }$E_{i}$, $D_{i}$\textit{\ is a non-empty compact
subset of }$X_{i}$\textit{\ and denote }$X=\prod\limits_{i\in I}X_{i},$ $%
D=\prod\limits_{i\in I}D_{i}$\textit{;}

ii)\textit{\ }$A_{i}:X\rightarrow 2^{D_{i}}$\textit{\ is upper
semicontinuous such that for each }$x\in X,$\textit{\ }$A_{i}(x)$\textit{\
is a non-empty closed convex subset of }$X_{i}$\textit{;}

iii)\textit{\ }$P_{i}:X\times X\rightarrow 2^{X_{i}}$\textit{\ is upper
semicontinuous such that for each }$x\in X,$\textit{\ }$P_{i}(x)$\textit{\
is a non-empty closed convex subset of }$X_{i}$\textit{;}

iv)\textit{\ the set }$W_{i}=\left\{ (x,y)\in X\times X:A_{i}\left( x\right)
\cap P_{i}(x,y)\neq \emptyset \right\} $\textit{\ is open;}

v)\textit{\ for each }$x\in W_{i},$\textit{\ }$x_{i}\notin P_{i}(x,y)$%
\textit{;}

\textit{Then, there exists }$(x^{\ast },y^{\ast })\in D\times D$ \textit{%
such that }$x_{i}^{\ast }\in A_{i}\left( x^{\ast }\right) ,$\textit{\ }$%
y^{\ast }\in F_{i}(x^{\ast })$ \textit{and }$A_{i}\left( x^{\ast }\right)
\cap P_{i}\left( x^{\ast },y^{\ast }\right) =\emptyset $\textit{\ for each }$%
i\in I.$

\textit{Proof. }According to the assumption i), $D$ is a non-empty compact
subset of $X$. For each $i\in I,$ the correspondence $G_{i}:X\times
X\rightarrow 2^{D_{i}}$ $\ $defined by $G_{i}\left( x,y\right) =A_{i}\left(
x\right) \cap P_{i}\left( x,y\right) $ for each $(x,y)\in X\times X$ has
closed convex values and it is upper semicontinuous. Let us define $%
T_{i}:X\times X\rightarrow 2^{D_{i}}$ by

$T_{i}\left( x,y\right) =\left\{ 
\begin{array}{c}
(A_{i}\left( x\right) \cap P_{i}(x,y))\times F_{i}(x),\text{ if }x\in W_{i};
\\ 
A_{i}\left( x\right) \times F_{i}(x),\text{ \ \ \ \ \ \ \ \ \ \ \ \ \ if }%
x\notin W_{i}.%
\end{array}%
\right. $

According to Lemma 3, $T_{i}$ is upper semicontinuous and has non-empty
closed convex values.

Let us define the correspondence $T:X\times X\rightarrow 2^{D}$, by $%
T(x,y)=\prod\limits_{i\in I}T_{i}(x,y)$ for each $x\in X.$ Then, $T$ is also
upper semicontinuous. Since each $T(x,y)$ is a non-empty closed convex
subset of the compact set $D,$ by Himmelberg's fixed point Theorem \cite{him}%
, there exists an $(x^{\ast },y^{\ast })\in D\times D$ such that $(x^{\ast
},y^{\ast })\in T\left( x^{\ast },y^{\ast }\right) .$ Then, $x_{i}^{\ast
}\in T_{i}\left( x^{\ast }\right) $ for each $i\in I.$ According to the
assumption v), $x_{i}^{\ast }\notin P_{i}(x^{\ast },y^{\ast })$ and
therefore, $x_{i}^{\ast }\in A_{i}\left( x^{\ast }\right) ,$ $y_{i}^{\ast
}\in F_{i}(x^{\ast })$ and $(A_{i}\cap P_{i})\left( x^{\ast }\right)
=\emptyset $ for each $i\in I.$

\begin{theorem}
\textit{Let }$\Gamma =\left\{ X_{i},A_{i},\text{ }F_{i},P_{i}\right\} _{i\in
I}$\textit{\ be a generalized abstract economy, where }$I$\textit{\ is any
index set, such that, for each }$i\in I$\textit{:}
\end{theorem}

i)\textit{\ }$X_{i}$\textit{\ is a non-empty convex subset of a Hausdorff
locally convex space }$E_{i}$, $D_{i}$\textit{\ is a non-empty compact
subset of }$X_{i}$\textit{\ and denote }$X=\prod\limits_{i\in I}X_{i},$ $%
D=\prod\limits_{i\in I}D_{i}$\textit{;}

ii)\textit{\ }$\overline{A_{i}^{V_{i}}},\overline{F_{i}^{V_{i}}}%
:X\rightarrow 2^{D_{i}}$\textit{\ are correspondences with non-empty convex
values for each open absolutely convex symmetric neighborhood }$V_{i}$%
\textit{\ of }$0$\textit{\ in }$E_{i}$;

iii)\textit{\ }$P_{i}:X\times X\rightarrow 2^{X_{i}}$\textit{\ is upper
semicontinuous with non-empty closed convex values;}

iv)\textit{\ the set }$W_{i}=\left\{ (x,y)\in X\times X:A_{i}\left( x\right)
\cap P_{i}(x,y)\neq \emptyset \right\} $\textit{\ is open;}

v)\textit{\ for each }$(x,y)\in W_{i},$\textit{\ }$x_{i}\notin P_{i}(x,y)$%
\textit{;}

\textit{Then there exists }$(x^{\ast },y^{\ast })\in D\times D$\textit{\
such that }$x_{i}^{\ast }\in \overline{A}_{i}\left( x^{\ast }\right) ,$ $%
y_{i}^{\ast }\in \overline{F}_{i}\left( x^{\ast }\right) $\textit{\ and }$%
A_{i}(x^{\ast })\cap P_{i}(x^{\ast },y^{\ast })=\emptyset $ for each $i\in
I. $

\textit{Proof.} According to assumption iv), $W_{i}$ is open in $X$ for each 
$i\in I.$

Let \ss $_{i}$\ be a basis of open absolutely convex symmetric neighborhoods
of $0$ in $E_{i}$ and let \ss =$\tprod\limits_{i\in I}$\ss $_{i}.$

For each $V=(V_{i})_{i\in I}\in \tprod\limits_{i\in I}$\ss $_{i},$ for each $%
i\in I,$ let's define $T_{i}^{V_{i}}:X\times X\rightarrow 2^{X_{i}}$ by

$T_{i}^{V_{i}}\left( x,y\right) :=\left\{ 
\begin{array}{c}
(\overline{A_{i}^{V_{i}}}\left( x\right) \cap P_{i}\left( x,y\right) )\times 
\overline{F_{i}^{V_{i}}}(x),\text{ if }(x,y)\in W_{i}, \\ 
\overline{A_{i}^{V_{i}}}\left( x\right) \times \overline{F_{i}^{V_{i}}}(x),%
\text{ \ \ \ \ \ \ \ \ \ \ \ if }(x,y)\notin W_{i}%
\end{array}%
\right. $ for each $(x,y)\in X\times X.$

According to assumption ii), each $T_{i}^{V_{i}}$ is upper semicontinuous
with non-empty closed convex values. Let us define $T^{V}:X\times
X\rightarrow 2^{D}$ by $T^{V}(x,y)=\tprod\limits_{i\in I}T_{i}^{V_{i}}(x,y)$
for each $(x,y)\in X\times X.$ The correspondence $T^{V}$ is upper
semicontinuous with non-empty closed convex values. Therefore, according to
Himmelberg's fixed point Theorem \cite{him}, there exists $(x_{V}^{\ast
},y_{V}^{\ast })=\tprod\limits_{i\in I}(x_{V}^{\ast },y_{V}^{\ast })_{i}\in
D\times D$ such that $(x_{V}^{\ast },y_{V}^{\ast })\in T^{V}(x_{V}^{\ast
},y_{V}^{\ast }).$ It follows that $(x_{V}^{\ast },,y_{V}^{\ast })_{i}\in 
\overline{T_{i}^{V_{i}}}(x_{V}^{\ast },,y_{V}^{\ast })$ for each $i\in I.$

For each $V=(V_{i})_{i\in I}\in $\ss $,$ let's define $Q_{V}=\cap _{i\in
I}\{(x,y)\in D\times D:$ $(x_{i},y_{i})\in T_{i}^{V_{i}}(x,y)\}.$

$Q_{V}$ is nonempty since $(x_{V}^{\ast },y_{V}^{\ast })\in Q_{V},$ then $%
Q_{V}$ is nonempty and closed.

We prove that the family $\{Q_{V}:V\in \text{\ss }\}$ has the finite
intersection property.

Let $\{V^{(1)},V^{(2)},...,V^{(n)}\}$ be any finite set of $\text{\ss }$ and
let $V^{(k)}=\underset{i\in I}{\tprod }V_{i}^{(k)}$, $k=1,...,n.$ For each $%
i\in I$, let $V_{i}=\underset{k=1}{\overset{n}{\cap }}V_{i}^{(k)}$, then $%
V_{i}\in \text{\ss }_{i};$ thus $V=\underset{i\in I}{\tprod }V_{i}\in 
\underset{i\in I}{\tprod }\text{\ss }_{i}.$ Clearly $Q_{V}\subseteq \underset%
{k=1}{\overset{n}{\cap }}Q_{V^{(k)}}$ so that $\underset{k=1}{\overset{n}{%
\cap }}Q_{V^{(k)}}\neq \emptyset .$

Since $D\times D$ is compact and the family $\{Q_{V}:V\in \text{\ss }\}$ has
the finite intersection property, we have that $\cap \{Q_{V}:V\in \text{\ss }%
\}\neq \emptyset .$ Take any $(x^{\ast },y^{\ast })\in \cap \{Q_{V}:V\in $%
\ss $\},$ then for each $V_{i}\in \text{\ss }_{i},$ $(x_{i}^{\ast
},y_{i}^{\ast })\in \overline{T_{i}^{V_{i}}}(x^{\ast },y^{\ast })$.
According to Lemma 2,\emph{\ }we have that\emph{\ } $(x_{i}^{\ast
},y_{i}^{\ast })\in \overline{T_{i}}(x^{\ast },y^{\ast }),$ for each $i\in
I. $

According to condition (5) we have that $x_{i}^{\ast }\in \overline{A}%
_{i}\left( x^{\ast }\right) ,$ $y_{i}^{\ast }\in \overline{F}_{i}\left(
x^{\ast }\right) $ and $(A_{i}\cap P_{i})(x^{\ast },y^{\ast })=\emptyset $
for each $i\in I.$

\section{Systems of vector quasi-equilibrium problems under upper
semicontinuity assumptions}

For each $i\in I,$ let $X_{i}$ be a non-empty subset of a topological vector
space $E_{i},$ $Y_{i}$ a topological vector space and let $%
X=\tprod\limits_{i\in I}X_{i}$ and $C_{i}\subset X_{i}$ a closed cone with
int$C\neq \emptyset .$ For each $i\in I,$ let $A_{i},$ $F_{i}:X\rightarrow
2^{X_{i}}$ and $f_{i}:X\times X\times X_{i}\rightarrow 2^{X_{i}}$ be
corespondences with non-empty values. We consider the following systems of
generalized vector quasi-equilibrium problems (in short, SGVQEP (I)):

Find $(x^{\ast },y^{\ast })\in X\times X$ such that for each $i\in I,$ $%
x_{i}^{\ast }\in \overline{A_{i}}(x^{\ast }),$ $y_{i}^{\ast }\in \overline{%
F_{i}}(x^{\ast })$ and $f_{i}(x^{\ast },y^{\ast },u_{i})\subseteq C_{i}$ for
each $u_{i}\in A_{i}(x^{\ast }).$

The next theorems establish the existence of the solutions for systems of
vector quasi-equil problem.

\begin{theorem}
For each $i\in I$ $(I$ finite)$,$ let $X_{i}$ be a non-empty compact convex
subset of a locally convex Hausdorff topological vector space $E_{i}$. Let $%
X=\tprod\limits_{i\in I}X_{i}$ \textit{be perfectly normal }and $C_{i}$ a
closed cone with int$C_{i}\neq \emptyset .$ Let $f_{i}:X\times X\times
X_{i}\rightarrow 2^{X_{i}}$ be a correspondence with non-empty values. For
each $i\in I,$ assume that:
\end{theorem}

\textit{i) }$F_{i}$\textit{, }$A_{i}:X\rightarrow 2^{X_{i}}$\textit{\ are
upper semicontinuous correspondences with non-empty closed convex values;}

\textit{ii) \ for all }$x,$\textit{\ }$y\in X,$\textit{\ }$%
f_{i}(x,y,x_{i})\subseteq C_{i};$

\textit{iii) }$f_{i}(\cdot ,\cdot ,\cdot )$\textit{\ is lower }$(-C_{i})-$%
\textit{semicontinuous;}

i\textit{v) for each }$x,y\in X,$\textit{\ the correspondence }$%
f_{i}(x,y,\cdot )$\textit{\ is }$C_{i}-$\textit{quasi-convex;}

\textit{v) }$U_{i}=\{(x,y)\in X\times X:$\textit{\ there exists }$u_{i}\in
A_{i}(x)$\textit{\ such that }$f_{i}(x,y,u_{i})\subseteq C_{i}\}$\textit{\
is open.}

\textit{Then, there exists a solution }$(x^{\ast },y^{\ast })\in X\times X$%
\textit{\ of (SGVQEP)(I).}

\textit{Proof. For each }$i\in I,$ let $P_{i}:X\times X\rightarrow 2^{X_{i}}$
be defined by

$P_{i}(x,y)=\{u_{i}\in X_{i}:f_{i}(x,y,u_{i})\varsubsetneq C_{i}\}$ for each 
$(x,y)\in X\times X.$

We will show that $P_{i}$ has an open graph and convex values.

We are proving firstly the convexity of $P_{i}(x_{0},y_{0}),$ where $%
(x_{0},y_{0})\in X\times X$ is arbitrary fixed$.$ Let us consider $%
u_{1},u_{2}\in P_{i}(x_{0},y_{0})$ and $\lambda \in \lbrack 0,1].$ Since $%
u_{1},u_{2}\in X_{i}$ and the set $X_{i}$ is convex, the convex combination $%
u=\lambda u_{1}+(1-\lambda )u_{2}\in X_{i}.$ Further, by using the property
of properly $C$-quasiconvexity of $f_{i}(x_{0},y_{0},\cdot )$, we can
assume, without loss of generality, that $f_{i}(x_{0},y_{0},u_{1})\subset
f_{i}(x_{0},y_{0},u)+C_{i}.$ We will prove that $u\in P_{i}(x_{0},y_{0}).$
If, by contrary, $u\notin P_{i}(x_{0},y_{0}),$ then, $f_{i}(x_{0},y_{0},u)%
\subseteq C_{i}$ and, consequently, $f_{i}(x_{0},y_{0},u_{1})\subset
f_{i}(x_{0},y_{0},u)+C_{i}\subseteq C_{i}+C_{i}\subseteq C_{i},$ which
contradicts $u_{1}\in P_{i}(x_{0},y_{0}).$ It remains that $u\in
P_{i}(x_{0},y_{0}).$ Therefore, $P_{i}(x_{0},y_{0})$ is a convex set.

According to assumption $ii)$, it follows that $x_{i}\notin P_{i}(x,y)$ for
each $(x,y)\in X\times X.$

The closedness of the $($Gr$P_{i})^{C}$ will be shown now. We consider the
net $\{(x_{\alpha },y_{\alpha },u_{\alpha }):\alpha \in \Lambda \}\subset ($%
Gr$P_{i})^{C}$ such that $(x_{\alpha },y_{\alpha },u_{\alpha })\rightarrow
(x_{0},y_{0},u_{0})\in X\times X\times X_{i}.$ Then, $u_{\alpha }\notin
P_{i}(x_{\alpha },y_{\alpha })$ for each $\alpha \in \Lambda $, i.e.$%
f_{i}(x_{\alpha },y_{\alpha },u_{\alpha })\subseteq C_{i}.$ We prove that $%
(x_{0},y_{0},u_{0})\in ($Gr$P_{i})^{C},$ that is $u_{0}\notin
P_{i}(x_{0},y_{0}).$ We use the lower $(-C_{i})-$ continuity of $F$ and we
conclude that, for each neighborhood $U$ of the origin in $X_{i},$ there
exists a neighbourhood $V(x_{0},y_{0},u_{0})$ of $(x_{0},y_{0},u_{0})$ such
that, $f_{i}(x_{0},y_{0},u_{0})\subset f_{i}(x_{\alpha },y_{\alpha
},u_{\alpha })+U+C_{i}$ for each $(x_{\alpha },y_{\alpha },u_{\alpha })\in
V(x_{0},y_{0},u_{0}).$ Then, for each $(x,y,u)\in V(x_{0},y_{0},u_{0}),$ $%
f_{i}(x_{0},y_{0},u_{0})\subset C_{i}+U+C_{i}\subset C_{i}+U.$ We will prove
that $f_{i}(x_{0},y_{0},u_{0})\subseteq C_{i}.$ If, by contrary, there
exists $a\in f_{i}(x_{0},y_{0},u_{0})$ and $a\notin C_{i},$ then, $0\notin
B:=C_{i}-a$ and $B$ is closed. Thus, $B$ is open and $0\in X_{i}\backslash
B. $ There exists an open symmetric neighborhood $U_{1}$ of the origin in $%
X_{i},$ such that $U_{1}\subset X_{i}\backslash B$ and $U_{1}\cap B$ is
closed. Therefore, $0\notin B+U_{1}$, i.e., $a\notin C_{i}+U_{1},$ which
contradicts $f(x_{0},y_{0},z_{0})\subset U_{1}+C_{i}.$ It follows that $%
f(x_{0},y_{0},z_{0})\subseteq C_{i}$ and then, $u_{0}\in
P_{i}(x_{0},y_{0})^{C}$ and $($Gr$P_{i})^{C}$ is closed and therefore, Gr$%
P_{i}$ is open and $P_{i}$ has open lower sections.

According to $vi)$, $U_{i}=\{(x,y)\in X\times X:$ $A_{i}(x)\cap
P_{i}(x,y)\neq \emptyset \}$ is open and according to $ii),$ $x_{i}\notin
P_{i}(x,y)$ for each $(x,y)\in X\times X.$

All the assumptions of Theorem 5 are fulfilled. Then, there exists $(x^{\ast
},y^{\ast })\in X\times X$ such that for each $i\in I,$ $A_{i}(x^{\ast
})\cap P_{i}(x^{\ast },y^{\ast })=\emptyset $, $x_{i}^{\ast }\in \overline{%
A_{i}}(x^{\ast })$ and\textit{\ }$y_{i}^{\ast }\in \overline{F_{i}}(x^{\ast
}).$ Consequently, there exist $x^{\ast },y^{\ast }\in X$ such that $%
x_{i}^{\ast }\in \overline{A_{i}}(x^{\ast })$,\textit{\ }$y_{i}^{\ast }\in 
\overline{F_{i}}(x^{\ast })$ and $f_{i}(x^{\ast },y^{\ast },u)\subseteq
C_{i} $ for each $u\in A_{i}(x^{\ast }).$

\begin{remark}
We note that Theorem 8 differs from Theorem 3.2.1 in \cite{lin3} in the
following way: the correspondences $A_{i}$ and $F_{i}$ are upper
semi-continuous and $f_{i}$ is lower $(-C_{i})-$ continuous for each $i\in I$%
.
\end{remark}

\begin{theorem}
For each $i\in I,$ let $X_{i}$\textit{\ be a non-empty convex subset of a
Hausdorff locally convex space }$E_{i}$, $D_{i}$\textit{\ a non-empty
compact subset of }$X_{i}$\textit{\ and denote }$X=\prod\limits_{i\in
I}X_{i},$ $D=\prod\limits_{i\in I}D_{i}$\textit{. }Let $f_{i}:X\times
X\times X_{i}\rightarrow 2^{X_{i}}$ be correspondence with non-empty values.
For each $i\in I,$ assume that:
\end{theorem}

\textit{i) }$F_{i}$\textit{, }$A_{i}:X\rightarrow 2^{X_{i}}$\textit{\ are
upper semicontinuous correspondences with non-empty closed convex values;}

\textit{ii) \ for each }$x,$\textit{\ }$y\in X,$\textit{\ }$%
f_{i}(x,y,x_{i})\subseteq $int$C_{i};$

\textit{iii) for each }$(x,y)\in X\times X,$\textit{\ }$f_{i}(x,y,\cdot )$%
\textit{\ is upper\ semicontinuous;}

\textit{iv) }$f_{i}(\cdot ,\cdot ,\cdot )$\textit{\ is upper (}$C_{i})-$%
\textit{\ semicontinuous;}

\textit{v) for each }$x,y\in X,$\textit{\ the correspondence }$%
f_{i}(x,y,\cdot )$\textit{\ is }$C_{i}-$\textit{quasi-convex;}

\textit{vi) }$U_{i}=\{(x,y)\in X\times X:$\textit{\ there exists }$u_{i}\in
A_{i}(x)$\textit{\ such that }$f(x,y,u_{i})\subseteq $int$C_{i}\}$\textit{\
is open.}

\textit{Then, there exists }$(x^{\ast },y^{\ast })\in X\times X$\textit{\
such that }$x_{i}^{\ast }\in \overline{A_{i}}(x^{\ast })$,\textit{\ }$%
y_{i}^{\ast }\in \overline{F_{i}}(x^{\ast })$ and $f_{i}(x^{\ast },y^{\ast
},u)\subseteq $int$C_{i}$ for each $u\in A_{i}(x^{\ast })$\textit{.}

\textit{Proof.} \textit{For each }$i\in I,$ let $P_{i}:X\times X\rightarrow
2^{X_{i}}$ be defined by

$P_{i}(x,y)=\{u_{i}\in X_{i}:f_{i}(x,y,u_{i})\varsubsetneq $int$C_{i}\}$ for
each $(x,y)\in X\times X.$

We will prove that $P_{i}$ has a closed graph and non-empty closed convex
values.

Let us fix $(x_{0},y_{0})\in X\times X$.

In order to prove the convexity of $P_{i}(x_{0},y_{0}),$ let us consider $%
u_{1},u_{2}\in P_{i}(x_{0},y_{0})$ and $\lambda \in \lbrack 0,1].$ Let $u$
be the convex combination $u=\lambda u_{1}+(1-\lambda )u_{2}\in X_{i}.$
Further, by using the property of properly C-quasiconvexity of $%
f_{i}(x_{0},y_{0},\cdot )$, we can assume, without loss of generality, that $%
f_{i}(x_{0},y_{0},u_{1})\subset f_{i}(x_{0},y_{0},u)+C_{i}.$ We will prove
that $u\in P_{i}(x_{0},y_{o}).$ If, by contrary, $u\notin
P_{i}(x_{0},y_{0}), $ $f_{i}(x_{0},y_{0},u)\subseteq $int$C_{i}$ and,
consequently, $f_{i}(x_{0},y_{0},u_{1})\subset
f_{i}(x_{0},y_{0},u)+C_{i}\subseteq $int$C_{i}+C_{i}\subseteq $int$C_{i},$
which contradicts $u_{1}\in P_{i}(x_{0},y_{0}).$ It remains that $u\in
P_{i}(x_{0},y_{0}).$ Therefore, $P_{i}(x_{0},y_{0})$ is a convex set.

Further, we will prove that $P_{i}(x_{0},y_{0})$ is closed$.$

Let us consider the net $\{u_{\alpha }:\alpha \in \Lambda \}\subseteq
P_{i}(x_{0},y_{0})$ such that $u_{\alpha }\rightarrow u_{0}.$ Then, $%
u_{\alpha }\in X_{i}$ and $f_{i}(x_{0},y_{0},u_{\alpha })\varsubsetneq $int$%
C_{i}$ for all $\alpha \in \Lambda .$ Since $X_{i}$ is a closed set, it
follows that $u_{0}\in X_{i}.$ We assume, by absurd, that $%
f_{i}(x_{0},y_{0},u_{0})\subseteq $int$C_{i}.$ Since $f_{i}(x_{0},y_{0},%
\cdot )$ is upper semicontinuous, then, $f_{i}(x_{0},y_{0},u_{\alpha
})\subset $int$C_{i}$ for $\alpha \geq \alpha _{0},$ $\alpha _{0}\in \Lambda
,$ which is a contradiction. Therefore, our assumption is false and $%
f_{i}(x_{0},y_{0},u_{0})\varsubsetneq $int$C_{i},$ i.e. $u_{0}\in
P_{i}(x_{0},y_{0})$ and $P_{i}(x_{0},y_{0})$ is a closed set.

Now, the closedness of $P_{i}$ will be shown. We consider the net $%
\{(x_{\alpha },y_{\alpha },u_{\alpha }):\alpha \in \Lambda \}\subset $Gr$%
P_{i}$ such that $(x_{\alpha },y_{\alpha },u_{\alpha })\rightarrow
(x_{0},y_{0},u_{0})\in X\times X\times X_{i}.$ Then, $u_{\alpha }\in
P_{i}(x_{\alpha },y_{\alpha })$ for each $\alpha \in \Lambda $ and we prove
that $(x_{0},y_{0},u_{0})\in $Gr$P_{i},$ that is $u_{0}\in
P_{i}(x_{0},y_{0}).$ If, by absurd, $u_{0}\notin P_{i}(x_{0},y_{0}),$ then, $%
f_{i}(x_{0},y_{0},u_{0})\subseteq $int$C_{i}.$ This relation implies that
there exists a neighbourhood $U_{0}$ of the origin in $Z$ such that $%
f_{i}(x_{0},y_{0},u_{0})+U_{0}\subset $int$C_{i}.$ Further, we use the upper 
$C_{i}-$continuity of $f_{i}$ and we conclude that there exists a
neighbourhood $V(x_{0},y_{0},u_{0})$ of $(x_{0},y_{0},u_{0})$ such that, $%
f_{i}(x,y,u)\subset f_{i}(x_{0},y_{0},u_{0})+U_{0}+C$ for each $(x,y,u)\in
V(x_{0},y_{0},u_{0}).$ Then, for each $(x,y,u)\in V(x_{0},y_{0},u_{0}),$ $%
f_{i}(x,y,u)\subset $int$C_{i}+C_{i}\subset $int$C_{i},$ which implies the
existence of $\alpha _{0}\in \Lambda $ such that for each $\alpha \geq
\alpha _{0},$ $f_{i}(x_{\alpha },y_{\alpha },u_{\alpha })\subset $int$C_{i}.$
The last relation contradicts $u_{\alpha }\in P_{i}(x_{\alpha },y_{\alpha
}). $ Consequently, the assumption that $u_{0}\notin P_{i}(x_{0},y_{0})$ is
false. Since $u_{0}\in P_{i}(x_{0},y_{0}),$ Gr$P_{i}$ is closed, and, since $%
X_{i}$ is compact, it follows that $P_{i}$ is upper semicontinuous.

According to vi), $U_{i}=(x,y)\in X\times X:$ $A_{i}(x)\cap P_{i}(x,y)\neq
\emptyset \}$ is open and according to $ii),$ $x_{i}\notin P_{i}(x,y)$ for
each $(x,y)\in X\times X.$

All the assumptions of Theorem 6 are fulfilled. Then, there exists $(x^{\ast
},y^{\ast })$ such that for each $i\in I,$ $A_{i}(x^{\ast })\cap
P_{i}(x^{\ast },y^{\ast })=\emptyset $, $x_{i}^{\ast }\in \overline{A_{i}}%
(x^{\ast })$ and\textit{\ }$y_{i}^{\ast }\in \overline{F_{i}}(x^{\ast }).$
Consequently, there exist $x^{\ast },y^{\ast }\in X$ such that $x_{i}^{\ast
}\in \overline{A_{i}}(x^{\ast })$,\textit{\ }$y_{i}^{\ast }\in \overline{%
F_{i}}(x^{\ast })$ and $f(x^{\ast },y^{\ast },u)\subseteq $int$C_{i}$ for
each $u\in A_{i}(x^{\ast }).$

\begin{theorem}
For each $i\in I,$ let $X_{i}$\textit{\ be a non-empty convex subset of a
Hausdorff locally convex space }$E_{i}$, $D_{i}$\textit{\ a non-empty
compact subset of }$X_{i}$\textit{\ and denote }$X=\prod\limits_{i\in
I}X_{i},$ $D=\prod\limits_{i\in I}D_{i}$\textit{. }Let $f_{i}:X\times
X\times X_{i}\rightarrow 2^{X_{i}}$ be a lower semicontinuous correspondence
with non-empty values. For each $i\in I,$ assume that:
\end{theorem}

\textit{i) }$\overline{F_{i}^{V_{i}}}$\textit{, }$\overline{A_{i}^{V_{i}}}%
:X\rightarrow 2^{X_{i}}$\textit{\ are correspondences with non-empty convex
values for each open absolutely convex symmetric neighborhood }$V_{i}$%
\textit{\ of }$0$\textit{\ in }$E_{i}$;

\textit{ii) \ for all }$x,$\textit{\ }$y\in X,$\textit{\ }$%
f_{i}(x,y,x_{i})\subseteq $int$C_{i};$

\textit{iii)} \textit{for each }$(x,y)\in X\times X,$\textit{\ }$%
f_{i}(x,y,\cdot )$\textit{\ is upper semicontinuous;}

\textit{iv) }$f_{i}(\cdot ,\cdot ,\cdot )$\textit{\ is upper (}$C_{i})-$%
\textit{\ semicontinuous;}

\textit{v) for each }$x,y\in X,$\textit{\ the correspondence }$%
f_{i}(x,y,\cdot )$\textit{\ is }$C_{i}-$\textit{quasi-convex;}

\textit{vi) }$U_{i}=\{(x,y)\in X\times X:$\textit{\ there exists }$u_{i}\in
A_{i}(x)$\textit{\ such that }$f(x,y,u_{i})\subseteq $int$C_{i}\}$\textit{\
is open.}

\textit{Then, there exists }$(x^{\ast },y^{\ast })\in X\times X$\textit{\
such that }$x_{i}^{\ast }\in \overline{A_{i}}(x^{\ast })$,\textit{\ }$%
y_{i}^{\ast }\in \overline{F_{i}}(x^{\ast })$ and $f_{i}(x^{\ast },y^{\ast
},u)\subseteq $int$C_{i}$ for each $u\in A_{i}(x^{\ast })$\textit{.}

\textit{Proof. }For each\textit{\ }$i\in I,$ let $P_{i}:X\times X\rightarrow
2^{X_{i}}$ be defined by

$P_{i}(x,y)=\{u_{i}\in X_{i}:f_{i}(x,y,u_{i})\varsubsetneq $int$C_{i}\}$ for
each $(x,y)\in X\times X.$

As in the proof of Theorem 9, we can show that $P_{i}$ is upper
semicontinuous with non-empty closed convex values.

According to vi), $U_{i}=(x,y)\in X\times X:$ $A_{i}(x)\cap P_{i}(x,y)\neq
\emptyset \}$ is open and according to $ii),$ $x_{i}\notin P_{i}(x,y)$ for
each $(x,y)\in X\times X.$

All the assumptions of Theorem 7 are fulfilled. Then, there exists $(x^{\ast
},y^{\ast })$ such that for each $i\in I,$ $A_{i}(x^{\ast })\cap
P_{i}(x^{\ast },y^{\ast })=\emptyset $, $x_{i}^{\ast }\in \overline{A_{i}}%
(x^{\ast })$ and\textit{\ }$y_{i}^{\ast }\in \overline{F_{i}}(x^{\ast }).$
Consequently, there exists $(x^{\ast },y^{\ast })\in X\times X$\ such that $%
x_{i}^{\ast }\in \overline{A_{i}}(x^{\ast })$,\ $y_{i}^{\ast }\in \overline{%
F_{i}}(x^{\ast })$ and $f_{i}(x^{\ast },y^{\ast },u)\subseteq $int$C_{i}$
for each $u\in A_{i}(x^{\ast })$\textit{.}

\section{Strong vector quasi-equilibrium problems}

Let us consider $E_{1},E_{2}$ and $Z$ be topological vector spaces, let $%
K\subset E_{1},$ $D\subset E_{2}$ be subsets and $C_{i}\subset Z$ a nonempty
closed convex cone. Let us also consider the correspondences $A:K\rightarrow
2^{K},$ $F:K\rightarrow 2^{D}$ and $f:K\times D\times K\rightarrow 2^{Z}.$

We will study the existence of the solutions for the following extension of
the generalized strong vector quasi-equilibrium problem (shortly, GSVQEP):
finding $x^{\ast }\in K$ and $y^{\ast }\in \overline{F}(x^{\ast })$ such
that $x^{\ast }\in \overline{A}(x^{\ast })$ and $f_{i}(x^{\ast },y^{\ast
},z)\subset C_{i},$ $\forall z\in \overline{A}(x^{\ast }),$ where the
correspondence $\overline{A}$ is defined by $\overline{A}(x)=\{y\in
Y:(x,y)\in $cl$_{X\times Y}$Gr$A\}.$ Note that cl$A(x)\subset \overline{A}%
(x) $ for each $x\in X.$

The element $x^{\ast }$ will be called a strong solution for the GSVQEP and
the set of all strong solutions for the GSVQEP will be denoted by $V_{A}(f).$

\subsection{Strong vector quasi-equilibrium problems without continuity
assumptions}

In this subsection, we prove that the set $V_{A}(f)$ is non-empty also in
the case when the correspondences do not satisfy continuity assumptions
defined explicitely. Instead, we will use some conditions concerning
generalized convexity, mainly, the weakly naturally quasi-concavity property
of correspondences. Other assumptions refer to correspondences with weakly
convex graphs. We present these notions below.

Let us denote $\Delta _{n-1}=\left\{ (\lambda _{1},\lambda _{2},...,\lambda
_{n})\in \mathbb{R}^{n}:\overset{n}{\underset{i=1}{\tsum }}\lambda _{i}=1%
\text{ and }\lambda _{i}\geqslant 0,i=1,2,...,n\right\} $ the standard
(n-1)-dimensional simplex in $\mathbb{R}^{n}.$

The correspondence\textit{\ }$F:X\rightarrow 2^{Y}$ is said to have \emph{%
weakly convex graph} \cite{ding2} if for each $n\in N$ and for each finite
set $\{x_{1},x_{2},...,x_{n}\}\subset X$, there exists $y_{i}\in F(x_{i})$, $%
(i=1,2,...,n)$ such that%
\begin{equation*}
(1.1)\ \ \ \ \ \text{co}(\{(x_{1},y_{1}),(x_{2},y_{2}),...,(x_{n},y_{n})\})%
\subset \text{Gr}(F)
\end{equation*}%
The relation (1.1) is equivalent to

\begin{center}
\begin{equation*}
(1.2)\ \ \ \ \overset{n}{\underset{i=1}{\tsum }}\lambda _{i}y_{i}\in F(%
\overset{n}{\underset{i=1}{\tsum }}\lambda _{i}x_{i})\ \ \ \ \ \ \ (\forall
(\lambda _{1},\lambda _{2},...,\lambda _{n})\in \Delta _{n-1}).
\end{equation*}
\end{center}

We introduced in \cite{pat} the weakly naturally quasi-concave
correspondences.

\begin{definition}
(see \cite{pat})Let $X$ be a non-empty convex subset of a topological vector
space $E$ \ and $Y$ a non-empty subset of a topological vector space $Z$.%
\textit{\ }The set-valued map\textit{\ }$F:X\rightarrow 2^{Y}$ is said to be 
\emph{weakly naturally quasi-concave (WNQ) }iff for each $n$ and for each
finite set $\{x_{1},x_{2},...,x_{n}\}\subset X$, there exists $y_{i}\in
F(x_{i})$, $(i\in \{1,...,n\}$ and $g=(g_{1},g_{2},...,g_{n}):\Delta
_{n-1}\rightarrow \Delta _{n-1}$ a mapping with $g_{i}$ continuous, $%
g_{i}(1)=1$ and $g_{i}(0)=0$ for each\textit{\ }$i\in \{1,2,...,n\}$, such
that $\overset{n}{\underset{i=1}{\tsum }}g_{i}(\lambda _{i})y_{i}\in F(%
\overset{n}{\underset{i=1}{\tsum }}\lambda _{i}x_{i})$ for every $(\lambda
_{1},\lambda _{2},...,\lambda _{n})\in \Delta _{n-1}.$
\end{definition}

\begin{example}
(see [20]) Let $F:[0,4]\rightarrow 2^{[-2,2]}$ be defined by
\end{example}

$F(x)=\left\{ 
\begin{array}{c}
\lbrack 0,2]\text{ if }x\in \lbrack 0,2); \\ 
\lbrack -2,0]\text{ \ \ if \ }x=2; \\ 
(0,2]\text{ if }x\in (2,4].%
\end{array}%
\right. $

$F$ is neither upper semicontinuous, nor lower semicontinuous in $2.$ $F$ is
weakly naturally quasi-concave.

The next theorem is our first result concerning the existence of the strong
solutions for the GSVQEP.

\begin{theorem}
Let $E_{1},E_{2}$,$Z$ be topological vector spaces, $K\subset E_{1}$ and $%
D\subset E_{2}$ be subsets. \textit{Let }$L$\textit{\ be a simplex in }$%
K\times D$\textit{\ and denote }$L_{K}=$pr$_{K}L.$ \textit{Let }$%
(A,F):L_{K}\rightarrow 2^{L}$\textit{\ be weakly naturally quasi-concave.
Let us suppose that, for each }$n\in N,$\textit{\ }$\lambda \in \Delta _{n-1}
$\textit{\ and }$x_{1},x_{2},...,x_{n}\in L_{k},$\textit{\ }$%
A(\tsum\limits_{i=1}^{n}\lambda _{i}x_{i})\subset
\tbigcap\limits_{i=1}^{n}A(x_{i}).$ \textit{Let }$f:L\times K\rightarrow
2^{K}$\textit{\ a corespondence such that the following assumptions are
fulfilled:}
\end{theorem}

\textit{i) }$\forall (x,y)\in L,$\textit{\ }$f(x,y,A(x))\subset C;$

\textit{ii) for each }$z\in K,$\textit{\ for each }$n\in N,$\textit{\ }$%
\lambda ,\lambda ^{\prime }\in \Delta _{n-1}$\textit{\ and }$%
(x_{1},y_{1}),(x_{2},y_{2}),...,(x_{n},y_{n})\in L,$\textit{\ }$%
f(\tsum\limits_{i=1}^{n}\lambda _{i}x_{i},\tsum\limits_{i=1}^{n}\lambda
_{i}^{\prime }y_{i},z)\subset \tbigcap\limits_{i=1}^{n}f(x_{i},y_{i},z).$

\textit{Then, }$V_{f}\neq \emptyset .$

\textit{Proof. }Let us define $P:L\rightarrow 2^{K},$ $M:L\rightarrow 2^{L},$
$M(x,y)=(P(x,y),F(x)),$ where $P(x,y)=\{u\in A(x):f(u,y,z)\subset C$ $%
\forall z\in A(x)\}.$ Firstly, the weakly naturally quasi-concavity of $M$
will be proved. Let us consider $n\in N,$ $%
(x_{1},y_{1}),(x_{2},y_{2}),...,(x_{n},y_{n})\in L.$ For each $i=1,2,...,n,$
there exists $(u_{i},v_{i})\in M(x_{i},y_{i}),$ that is $u_{i}\in A(x_{i})$
and $f(u_{i},y_{i},z)\subset C$ $\forall z\in A(x_{i})$ and $v_{i}\in
F(x_{i}).$ Let $\lambda \in \Delta _{n-1}$ such that $\overset{n}{\underset{%
i=1}{\tsum }}\lambda _{i}(x_{i},y_{i})\in L$ and let us denote $(u,v)=%
\overset{n}{\underset{i=1}{\tsum }}\lambda _{i}(u_{i},v_{i}).$ Since $(A,F)$
is weakly naturally quasiconcave, it follows that there exists $%
g=(g_{1},g_{2},...,g_{n}):\Delta _{n-1}\rightarrow \Delta _{n-1}$ a function
depending on $(x_{1},y_{1}),(x_{2},y_{2}),...,(x_{n},y_{n})$ with $g_{i}$
continuous, $g_{i}(1)=1,$ $g_{i}(0)=0$ for each\textit{\ }$i=1,2,...n$, such
that for every $\lambda =(\lambda _{1},\lambda _{2},...,\lambda _{n})\in
\Delta _{n-1}$, $\overset{n}{\underset{i=1}{\tsum }}g_{i}(\lambda
_{i})(u_{i},v_{i})\in (A,F)(\overset{n}{\underset{i=1}{\tsum }}\lambda
_{i}x_{i}).$ According to assumption ii), $f(\overset{n}{\underset{i=1}{%
\tsum }}g_{i}(\lambda _{i})u_{i},\overset{n}{\underset{i=1}{\tsum }}\lambda
_{i}y_{i},z)\subset f(u_{i},y_{i},z)$

\noindent $\subset C.$ We have that $f(\overset{n}{\underset{i=1}{\tsum }}%
g_{i}(\lambda _{i})u_{i},\overset{n}{\underset{i=1}{\tsum }}\lambda
_{i}y_{i},z)\subset C$ $\forall z\in A(\overset{n}{\underset{i=1}{\tsum }}%
\lambda _{i}x_{i})$ and then, $f(\overset{n}{\underset{i=1}{\tsum }}%
g_{i}(\lambda _{i})u_{i},\overset{n}{\underset{i=1}{\tsum }}\lambda
_{i}y_{i},z)\subset C$ $\forall z\in A(x_{i}),$ because $A(\overset{n}{%
\underset{i=1}{\tsum }}\lambda _{i}x_{i})\subset A(x_{i})$ for each $i\in
\{1,2,...,n\}.$ Therefore, $P$ is weakly naturally quasi-concave.

Therefore, $\overset{n}{\underset{i=1}{\tsum }}g_{i}(\lambda _{i})u_{i}\in A(%
\overset{n}{\underset{i=1}{\tsum }}\lambda _{i}x_{i}),$ $f(\overset{n}{%
\underset{i=1}{\tsum }}g_{i}(\lambda _{i})u_{i},\overset{n}{\underset{i=1}{%
\tsum }}\lambda _{i}y_{i},z)\subset C$ $\forall z\in A(\overset{n}{\underset{%
i=1}{\tsum }}\lambda _{i}x_{i})$ and $\overset{n}{\underset{i=1}{\tsum }}%
g_{i}(\lambda _{i})v_{i}\in F(\overset{n}{\underset{i=1}{\tsum }}\lambda
_{i}x_{i}),$ that is $\overset{n}{\underset{i=1}{\tsum }}g_{i}(\lambda
_{i})u_{i}\in P(\overset{n}{\underset{i=1}{\tsum }}\lambda _{i}x_{i},\overset%
{n}{\underset{i=1}{\tsum }}\lambda _{i}y_{i})$ and $\overset{n}{\underset{i=1%
}{\tsum }}g_{i}(\lambda _{i})v_{i}\in F(\overset{n}{\underset{i=1}{\tsum }}%
\lambda _{i}x_{i}).$ Consequently, $\overset{n}{\underset{i=1}{\tsum }}%
g_{i}(\lambda _{i})(u_{i},v_{i})\in M(\overset{n}{\underset{i=1}{\tsum }}%
\lambda _{i}(x_{i},y_{i})$, hence, $M$ is weakly naturally quasi-concave.
Further, we will prove that $M$\textit{\ }has a continuous selection on%
\textit{\ }$L$\textit{. }$L$ is a simplex, and let us suppose that it is the
convex hull of the affinely independent set $%
\{(a_{1},b_{1}),(a_{2},b_{2}),...,(a_{n},b_{n})\}.$ There exist unique
continuous functions $\lambda _{i}:L\rightarrow \mathbb{R},$ $i=1,2,...,n$
such that, for each $(x,y)\in L,$ we have $(\lambda _{1}(x,y),\lambda
_{2}(x,y),...,\lambda _{n}(x,y))\in \Delta _{n-1}$ and $(x,y)=\overset{n}{%
\underset{i=1}{\tsum }}\lambda _{i}(x,y)(a_{i},b_{i}).$ Let us define $%
h:L\rightarrow L$ by $h(a_{i},b_{i})=(c_{i},d_{i})$ $(i=1,...,n)$ and $h(%
\overset{n}{\underset{i=1}{\tsum }}\lambda _{i}(a_{i},b_{i}))=\overset{n}{%
\underset{i=1}{\tsum }}\lambda _{i}(c_{i},d_{i})\in M(x,y).$ We show that $h$
is continuous. Let $(x_{m},y_{m})_{m\in N}$ be a sequence which converges to 
$(x_{0},y_{0})\in L,$ where $(x_{m},y_{m})=\overset{n}{\underset{i=1}{\tsum }%
}\lambda _{i}(x_{m},y_{m})(a_{i},b_{i})$ and $(x_{0},y_{0})=$ $\overset{n}{%
\underset{i=1}{\tsum }}\lambda _{i}(x_{0})(a_{i},b_{i}).$ By the continuity
of $\lambda _{i},$ it follows that, for each $i=1,2,...,n$, $\lambda
_{i}(x_{m},y_{m})\rightarrow \lambda _{i}(x_{0},y_{0})$ as $m\rightarrow
\infty .$ Hence, $h(x_{m},y_{m})\rightarrow h(x_{0},y_{0})$ as $m\rightarrow
\infty ,$ i.e. $h$ is continuous.

We proved that $M$ has a continuous selection on $B.$ According to Brouwer
fixed point Theorem, $h$ has a fixed point $(x^{\ast },y^{\ast })\in L,$ $%
h(x^{\ast },y^{\ast })=(x^{\ast },y^{\ast }).$ Then, $(x^{\ast },y^{\ast
})\in M(x^{\ast },y^{\ast }).$ Therefore, $x^{\ast }\in P(x^{\ast },y^{\ast
})$ and $y^{\ast }\in F(x^{\ast }),$ which implies that there exist $x^{\ast
}\in K$ and $y^{\ast }\in F(x^{\ast })$ such that $x^{\ast }\in A(x^{\ast })$
and $f(x^{\ast },y^{\ast },x)\subset C$, $\forall x\in A(x^{\ast }),$ i.e. $%
x^{\ast }\in V_{A}(f)$ and then, $V_{A}(f)$ is non-empty.\medskip

We will prove a similar result in the case of biconvexity. The biconvex sets
were introduced by Aumann \cite{aum}. For the reader's convenience, we
present below the most important notions concerning biconvexity.

Let $X\subset E_{1}$ and $Y\subset E_{2}$ be two nonempty convex sets, $%
E_{1},E_{2}$ be topological vector spaces and let $B\subset X\times Y.$

The set $B\subset X\times Y$ is called a \textit{biconvex set} on $X\times Y$
if the section $B_{x}=\left\{ y\in Y:(x,y)\in B\right\} $ is convex for
every $x\in X$ and the section $B_{y}=\left\{ x\in X:(x,y)\in B\right\} $ is
convex for every $y\in Y.$ Let $(x_{i},y_{i})\in X\times Y$ for $i=1,2,...n.$
A convex combination $(x,y)=\tsum\limits_{i=1}^{n}\lambda _{i}(x_{i},y_{i})$%
, (with $\tsum\limits_{i=1}^{n}\lambda _{i}=1,$ $\lambda _{i}\geq 0$ $%
i=1,2,...,n$) is called \textit{biconvex combination} if $%
x_{1}=x_{2}=...=x_{n}=x$ or $y_{1}=y_{2}=...=y_{n}=y.$ Let $D\subseteq
X\times Y$ be a given set. The set $H:=\tbigcap \{D_{I}:D\subseteq D_{I},$ $%
D_{I}$ is biconvex\} is called \textit{biconvex hull of }$D$ and is denoted
biconv$(D).\medskip $

\begin{theorem}
(Aumann and Hart \cite{aum}). A set $B\subseteq X\times Y$ is biconvex if
and only if $B$ contains all biconvex combinations of its elements.
\end{theorem}

\begin{theorem}
(Aumann and Hart \cite{aum}). The \textit{biconvex hull of a set} $P$ is
biconvex. Furthermore, it is the smallest biconvex set (in the sens of set
inclusion), which contains $P.$
\end{theorem}

\begin{lemma}
(Gorski, Pfeuffer and Klamroth \cite{gor}). Let $D\subseteq X\times Y$ be a
given set. Then biconv$(D)\subseteq $conv$(D).$
\end{lemma}

Now we introduce the following definition.

\begin{definition}
Let $B\subset X\times Y$ be a biconvex set, $Z$ a nonempty convex subset of
a topological vector space $F$ and $T:B\rightarrow 2^{Z}$ a correspondence$.$
$T$ is called \textit{weakly biconvex} if for each finite set $%
\{(x_{1},y_{1}),(x_{2},y_{2}),...,(x_{n},y_{n})\}\subset B$, there exists $%
z_{i}\in T(x_{i},y_{i})$, $(i=1,2,...,n)$ such that for every biconvex
combination $(x,y)=\tsum\limits_{i=1}^{n}\lambda _{i}(x_{i},y_{i})\in B$
(with $\tsum\limits_{i=1}^{n}\lambda _{i}=1,$ $\lambda _{i}\geq 0$ $%
i=1,2,...,n$)$,$ $\overset{n}{\underset{i=1}{\tsum }}\lambda _{i}z_{i}\in T(%
\overset{n}{\underset{i=1}{\tsum }}\lambda _{i}(x_{i},y_{i})).$
\end{definition}

\begin{theorem}
Let $E_{1},E_{2},Z$ be topological vector spaces and $K\subset E_{1},$ $%
D\subset E_{2}$ be subsets.
\end{theorem}

\textit{Let }$B$\textit{\ be the biconvex hull of }$%
\{(a_{1},b_{1}),(a_{2},b_{2}),...,(a_{n},b_{n})\}\subset K\times D$\textit{\
(a biconvex subset of }$K\times D)$\textit{\ and denote }$B_{K}=pr_{K}B.$

\textit{Let }$(A,F):B_{K}\rightarrow 2^{B},$\textit{\ such that and }$%
A:B_{K}\rightarrow 2^{K}$\textit{\ and }$F:B_{K}\rightarrow 2^{D}$\textit{\
have weakly convex graphs and for each }$n\in N,$\textit{\ }$\lambda \in
\Delta _{n-1}$\textit{\ and }$x_{1},x_{2},...,x_{n}\in B_{k},$\textit{\ }$%
A(\tsum\limits_{i=1}^{n}\lambda _{i}x_{i})\subset
\tbigcap\limits_{i=1}^{n}A(x_{i}).$

\textit{Let }$f:B\times K\rightarrow 2^{K}$\textit{\ such that:}

\textit{a) }$\forall (x,y)\in B,$\textit{\ }$f(x,y,A(x))\subset C$

\textit{b) for each }$z\in K,$\textit{\ for each }$n\in N,$\textit{\ }$%
\lambda \in \Delta _{n-1}$\textit{\ and }$%
(x_{1},y_{1}),(x_{2},y_{2}),...,(x_{n},y_{n})\in B,$\textit{\ }$%
f(\tsum\limits_{i=1}^{n}\lambda _{i}x_{i},\tsum\limits_{i=1}^{n}\lambda
_{i}y_{i},z)\subset \tbigcap\limits_{i=1}^{n}f(x_{i},y_{i},z).$

\textit{Then, }$V_{f}\neq \emptyset .$

\textit{Proof.} Let us define $P:B\rightarrow 2^{K}$ by $P(x,y)=\{u\in
A(x):f(u,y,z)\subset C$ $\forall z\in A(x)\}.$

We will prove that $P$ is weakly biconvex. Let $n\in N$ and $%
(x_{1},y_{1}),(x_{2},y_{2}),...,(x_{n},y_{n})\in B.$

For each $i=1,2,...,n,$ there exists $u_{i}\in P_{i}(x_{i},y_{i}),$ that is $%
u_{i}\in A(x_{i})$ and $f(u_{i},y_{i},z)\subset C$ $\forall z\in A(x_{i}).$
Let $\lambda \in \Delta _{n-1}$ such that $\tsum\limits_{i=1}^{n}\lambda
_{i}(x_{i},y_{i})\in B$ and let us denote $u=\tsum\limits_{i=1}^{n}\lambda
_{i}u_{i}.$ Since $A$ has a weakly convex graph, it follows that $u\in
A(\tsum\limits_{i=1}^{n}\lambda _{i}x_{i}).$

$f(\tsum\limits_{i=1}^{n}\lambda _{i}u_{i},\tsum\limits_{i=1}^{n}\lambda
_{i}y_{i},z)\subset f(u_{i},y_{i},z)\subset C$ $\forall z\in
A(\tsum\limits_{i=1}^{n}\lambda _{i}x_{i}).$ Since $A(\tsum\limits_{i=1}^{n}%
\lambda _{i}x_{i})\subset A(x_{i})$ for each $i\in \{1,2,...,n\},$ we have
that $f(\tsum\limits_{i=1}^{n}\lambda
_{i}u_{i},\tsum\limits_{i=1}^{n}\lambda _{i}y_{i},z)\subset C$ $\forall z\in
A(x_{i})$ $\forall z\in A(x_{i})$, $i=1,2,...,n.$ Therefore, $P$ is weakly
biconvex.

Let $M:B\rightarrow 2^{B},$ $M(x,y)=(P(x,y),F(x))$ $\forall (x,y)\in B.$ We
will prove that $M$ is weakly biconvex. Let $n\in N$ and $%
(x_{1},y_{1}),(x_{2},y_{2}),...,(x_{n},y_{n})\in B.$ Since $P$ is weakly
biconvex, for each $i=1,2,...,n,$ there exists $u_{i}\in P(x_{i},y_{i})$
such that for each $\lambda \in \Delta _{n-1},$ $\tsum\limits_{i=1}^{n}%
\lambda _{i}u_{i}\in P(\tsum\limits_{i=1}^{n}\lambda _{i}(x_{i},y_{i})).$
Since $F$ has a weakly convex graph, for each $i=1,2,...,n,$ there exists $%
v_{i}\in F(x_{i},y_{i})$ such that for each $\lambda \in \Delta _{n-1},$ $%
\tsum\limits_{i=1}^{n}\lambda _{i}v_{i}\in F(\tsum\limits_{i=1}^{n}\lambda
_{i}x_{i}).$ It follows that for each $i=1,2,...,n,$ there exists $%
(u_{i},v_{i})\in M(x_{i},y_{i})$ such that for each $\lambda \in \Delta
_{n-1},$ $\tsum\limits_{i=1}^{n}\lambda _{i}(u_{i},v_{i})\in
M(\tsum\limits_{i=1}^{n}\lambda _{i}(x_{i},y_{i})).$

Further, it will be shown that $M$\textit{\ }has a continuous selection on%
\textit{\ }$B$\textit{. }Since $B$ is biconvex hull of $%
(a_{1},b_{1}),...,(a_{n},b_{n}),$ there exist unique continuous functions $%
\lambda _{i}:K\rightarrow \mathbb{R},$ $i=1,2,...,n$ such that for each $%
(x,y)\in B,$ we have $(\lambda _{1}(x,y),\lambda _{2}(x,y),...,\lambda
_{n}(x,y))\in \Delta _{n-1}$ and $(x,y)=\overset{n}{\underset{i=1}{\tsum }}%
\lambda _{i}(x,y)(a_{i},b_{i}).$

Define $h:B\rightarrow B$ by $h(a_{i},b_{i})=(c_{i},d_{i})$ $(i=1,...,n)$
and $h(\overset{n}{\underset{i=1}{\tsum }}\lambda _{i}(a_{i},b_{i}))=\overset%
{n}{\underset{i=1}{\tsum }}\lambda _{i}(c_{i},d_{i})\in M(x,y).$ We show
that $h$ is continuous. Let $(x_{m},y_{m})_{m\in N}$ be a sequence which
converges to $(x_{0},y_{0})\in B,$ where $(x_{m},y_{m})=\overset{n}{\underset%
{i=1}{\tsum }}\lambda _{i}(x_{m},y_{m})(a_{i},b_{i})$ implies $%
a_{1}=a_{2}=...=a_{n}=a$ or $b_{1}=b_{2}=...=b_{n}=b$ and $(x_{0},y_{0})=$ $%
\overset{n}{\underset{i=1}{\tsum }}\lambda _{i}(x_{0})(a_{i},b_{i})$ with $%
a_{1}=a_{2}=...=a_{n}=a$ or $b_{1}=b_{2}=...=b_{n}=b.$ By the continuity of $%
\lambda _{i},$ it follows that for each $i=1,2,...,n$, $\lambda
_{i}(x_{m},y_{m})\rightarrow \lambda _{i}(x_{0},y_{0})$ as $m\rightarrow
\infty .$ Hence $h(x_{m},y_{m})\rightarrow h(x_{0},y_{0})$ as $m\rightarrow
\infty ,$ i.e. $h$ is continuous.

We proved that $M$ has a continuous selection on $B.$ According to Brouwer
fixed point Theorem$,$ $h$ has a fixed point $(x^{\ast },y^{\ast })\in B,$
i.e. $h(x^{\ast },y^{\ast })=(x^{\ast },y^{\ast }).$ It follows that $%
(x^{\ast },y^{\ast })\in M(x^{\ast },y^{\ast }).$ Therefore, $x^{\ast }\in
P(x^{\ast },y^{\ast })$ and $y^{\ast }\in F(x^{\ast }),$ which implies that
there exists $x^{\ast }\in K$ and $y^{\ast }\in F(x^{\ast })$ such that $%
x^{\ast }\in A(x^{\ast })$ and $f(x^{\ast },y^{\ast },x)\subset C$, $\forall
x\in A(x^{\ast }),$ i.e. $x^{\ast }\in V_{A}(f)$ and then, $V_{A}(f)$ is
non-empty.

\subsection{Strong vector quasi-equilibrium problems under weak continuity
assumptions}

Now, our last result is stated. The technique of proof is an approximation
one.

The next example shows that our condition over $A$ can be fulfilled by
correspondences which are not lower semicontinuous. Let us recall that, if $%
A:X\rightarrow 2^{Y}$ is a correspondence and $D,V\subset Y,$ then $%
A^{V}:X\rightarrow 2^{Y}$ is defined by $A^{V}(x)=(A(x)+V)\cap D$, $\forall
x\in X.$

\begin{example}
Let $\ A:(0,2)\rightarrow 2^{[1,4]}$ be the correspondence defined by
\end{example}

$A(x)=\left\{ 
\begin{array}{c}
\lbrack 2-x,2],\text{ if }x\in (0,1); \\ 
\{4\}\text{ \ \ \ \ \ \ if \ \ \ \ \ \ \ }x=1; \\ 
\lbrack 1,2]\text{ \ \ \ if \ \ \ }x\in (1,2).%
\end{array}%
\right. $ $A$ is not lower semicontinuous on $(0,2).$

Let $D=[1,2].$ For each $V=(-\varepsilon ,\varepsilon )$ with $\varepsilon
>0,$ the correspondence $\overline{A^{V}}$ is lower semicontinuous and $%
\overline{A^{V}}$ has nonempty convex values.

\begin{theorem}
Let $E_{1},E_{2},Z$ be Hausdorff locally convex topological vector spaces, $%
K\subset E_{1}$ and $D\subset E_{2}$ be non-empty convex compact subsets and 
$C$ be a non-empty closed convex cone. \textit{Let }$A:K\rightarrow 2^{K}$%
\textit{\ be a correspondence such that }$\overline{A}$ and $\overline{%
A^{V_{1}}}$\textit{\ are lower semi-continuous with non-empty convex values
for each absolutely convex symmetric neighbourhood }$V_{1}$\textit{\ of }$0$%
\textit{\ in }$X.$ \textit{Let }$F:K\rightarrow 2^{D}$\textit{\ be such that 
}$\overline{F}$ and $\overline{F^{V_{2}}}$\textit{\ has non-empty convex
values for each absolutely convex symmetric neighbourhood }$V_{2}$\textit{\
of }$0$\textit{\ in }$Y.$ \textit{Let }$f:K\times D\times K\rightarrow 2^{Z}$%
\textit{\ such that the following assumptions are satisfied:}
\end{theorem}

\textit{i) for all }$(x,y)\in K\times D,$\textit{\ }$f(x,y,\overline{A}%
(x))\subset C$\textit{\ and }$f(x,y,\overline{A^{V_{1}}}(x))\subset C$%
\textit{\ for each absolutely convex symmetric neighbourhood }$V_{1}$\textit{%
\ of }$0$\textit{\ in }$X;$

\textit{ii) for all }$(y,z)\in D\times K,$\textit{\ }$f(\cdot ,y,z)$\textit{%
\ is properly }$C-$\textit{quasiconvex;}

\textit{iii) }$f(\cdot ,\cdot ,\cdot )$\textit{\ is upper }$C-$\textit{%
continuous;}

\textit{iv) for all }$y\in D,$\textit{\ }$f(\cdot ,y,\cdot )$\textit{\ is
lower }$(-C)$\textit{-continuous.}

\textit{Then, }$V_{A}(f)\neq \emptyset .$

\textit{Proof. }Let us define $P:K\times D\rightarrow 2^{K}$ and $M:K\times
D\rightarrow 2^{K\times D}$ by

$P(x,y)=\{u\in A(x):f(u,y,z)\subset C,$ $\forall z\in A(x)\}$ $\forall
(x,y)\in K\times D$ and

$M(x,y)=(P(x,y),F(x))$ $\forall (x,y)\in K\times D.$

Let $V_{1}$ be an open absolutely convex symmetric neighbourhood of $0$ in $%
X.$

Firstly, we will show that $\overline{P^{V_{1}}}$ is an upper semicontinuous
correspondence with non-empty closed convex values. The elements of $%
V_{A}(f) $ will be obtained as a consequence of the existence of the fixed
points for $\overline{M}.$

Let $P^{V_{1}}$ be defined by $P^{V_{1}}(x,y)=\{u\in
A^{V_{1}}(x):f(u,y,z)\subset C,$ $\forall z\in A^{V_{1}}(x)\}$ $\forall
(x,y)\in K\times D.$

We will show that $\overline{P^{V_{1}}}$ is defined by $\overline{P^{V_{1}}}%
(x,y)=\{u\in \overline{A^{V_{1}}}(x):f(u,y,z)\subset C,$ $\forall z\in 
\overline{A^{V_{1}}}(x)\}$ $\forall (x,y)\in K\times D.$

In order to do this fact, the closedness of $\overline{P^{V_{1}}}$ will be
shown firstly. We consider the net $\{(x_{\alpha },y_{\alpha },u_{\alpha
}):\alpha \in \Lambda \}\subset $Gr$\overline{P^{V_{1}}}$ such that $%
(x_{\alpha },y_{\alpha },u_{\alpha })\rightarrow (x_{0},y_{0},u_{0})\in
K\times D\times K.$ Then, $u_{\alpha }\in \overline{P^{V_{1}}}(x_{\alpha
},y_{\alpha })$ for each $\alpha \in \Lambda $ and we prove that $%
(x_{0},y_{0},u_{0})\in $Gr$\overline{P^{V_{1}}},$ that is $u_{0}\in 
\overline{P^{V_{1}}}(x_{0},y_{0}).$ Since $\overline{A^{V_{1}}}$ is upper
semicontinuous and $u_{\alpha }\in \overline{A^{V_{1}}}(x_{\alpha }),$ then $%
u_{0}\in \overline{A^{V_{1}}}(x_{0}).$ If, by absurd, $u_{0}\notin \overline{%
P^{V_{1}}}(x_{0},y_{0}),$ there exists $z_{0}\in \overline{A^{V_{1}}}(x_{0})$
such that $f(u_{0},y_{0},z_{0})\varsubsetneq C.$ This relation implies that
there exists a neighbourhood $U_{0}$ of the origin in $Z$ such that $%
f(u_{0},y_{0},z_{0})+U_{0}\varsubsetneq C.$ Further, we use the upper $C-$
continuity of $f$ and we conclude that there exists a neighbourhood $%
V(u_{0},y_{0},z_{0})$ of $(u_{0},y_{0},z_{0})$ such that, $f(u,y,z)\subset
f(u_{0},y_{0},z_{0})+U_{0}+C$ for each $(u,y,z)\in V(u_{0},y_{0},z_{0}).$
Then, for each $(u,y,z)\in V(u_{0},y_{0},z_{0}),$ $f(u,y,z)\varsubsetneq
C+C\subset C,$ which implies the existence of $\alpha _{0}\in \Lambda $ such
that for each $\alpha \geq \alpha _{0},$ $f(u_{\alpha },y_{\alpha
},z_{\alpha })\varsubsetneq C.$ The last relation contradicts $u_{\alpha
}\in \overline{P^{V_{1}}}(x_{\alpha },y_{\alpha }).$ Consequently, the
assumption that $u_{0}\notin \overline{P^{V_{1}}}(x_{0},y_{0})$ is false.
Since $u_{0}\in \overline{P^{V_{1}}}(x_{0},y_{0}),$ Gr$\overline{P^{V_{1}}}$
is closed, and, since $K$ is compact, it follows that $\overline{P^{V_{1}}}$
is upper semicontinuous.

$\overline{A^{V_{1}}}(x)$ is non-empty for each $x\in K$ and, according to
the assumption i), it follows that $\overline{P^{V_{1}}}(x,y)$ is non-empty.

Let us fix $(x_{0},y_{0})\in K\times D$. We will prove secondly that $%
\overline{P^{V_{1}}}(x_{0},y_{0})$ is closed$.$ Let us consider the net $%
\{u_{\alpha }:\alpha \in \Lambda \}\subseteq \overline{P^{V_{1}}}%
(x_{0},y_{0})$ such that $u_{\alpha }\rightarrow u_{0}.$ Then, $u_{\alpha
}\in \overline{A^{V_{1}}}(x_{0})$ and $f(u_{\alpha },y_{0},z)\subset C$ for
all $z\in \overline{A^{V_{1}}}(u_{\alpha }).$ Since $\overline{A^{V_{1}}}%
(x_{0})$ is a closed set, it follows that $u_{0}\in \overline{A^{V_{1}}}%
(x_{0}).$ By using the lower semicontinuity of $\overline{A^{V_{1}}},$ we
obtain that, for any $z_{0}\in \overline{A^{V_{1}}}(u_{0})$ and $\{u_{\alpha
}\}\rightarrow u_{0},$ there exists a net $\{z_{\alpha }\}$ such that $%
z_{\alpha }\in \overline{A^{V_{1}}}(u_{\alpha })$ and $z_{\alpha
}\rightarrow z_{0}.$ It follows that, for each $\alpha ,$ there exists $%
z_{\alpha }\in \overline{A^{V_{1}}}(u_{\alpha })$ such that $f(u_{\alpha
},y_{0},z_{\alpha })\subset C.$ Since $f(\cdot ,y,\cdot )$ is lower $(-C)$%
-continuous, for each neighbourhood $U$ of the origin in $Z,$ there exists a
subnet $\{u_{\beta },z_{\beta }\}$ of $\{u_{\alpha },z_{\alpha }\}$ such
that $f(u_{0},y_{0},z_{0})\subset f(u_{\alpha },y_{0},z_{\alpha })+U+C.$
Consequently, $f(u_{0},y_{0},z_{0})\subset U+C.$ Further, we prove that $%
f(u_{0},y_{0},z_{0})\subset C.$ If, by contrary, there exists $a\in
f(u_{0},y_{0},z_{0})$ and $a\notin C,$ then, $0\notin B:=C-a$ and $B$ is
closed. Thus, $Z\backslash B$ is open and $0\in Z\backslash B.$ There exists
an open symmetric neighbourhood $U_{1}$ of the origin in $Z,$ such that $%
U_{1}\subset Z\backslash B$ and $U_{1}\cap B=\emptyset .$ Therefore, $%
0\notin B+U_{1},$ i.e., $a\notin C+U_{1},$ which contradicts $%
f(u_{0},y_{0},z_{0})\subset U_{1}+C$. It follows that $f(u_{0},y_{0},z_{0})%
\subset C$ and then, $u_{0}\in \overline{P^{V_{1}}}(x_{0},y_{0})$ and $%
\overline{P^{V_{1}}}(x_{0},y_{0})$ is closed.

Therefore, $\overline{P^{V_{1}}}(x,y)=\{u\in \overline{A^{V_{1}}}%
(x):f(u,y,z)\subset C,$ $\forall z\in \overline{A^{V_{1}}}(x)\}$ $\forall
(x,y)\in K\times D.$

Now, we are proving the convexity of $\overline{P^{V_{1}}}(x_{0},y_{0}),$
where $(x_{0},y_{0})\in X\times X$ is arbitrary fixed$.$ Let us consider $%
u_{1},u_{2}\in \overline{P^{V_{1}}}(x_{0},y_{0})$ and $\lambda \in \lbrack
0,1].$ Since $u_{1},u_{2}\in \overline{A^{V_{1}}}(x_{0})$ and the set $%
\overline{A^{V_{1}}}(x_{0})$ is convex, the convex combination $u=\lambda
u_{1}+(1-\lambda )u_{2}\in \overline{A^{V}}(x_{0}).$ Further, by using the
property of properly C-quasiconvexity of $f(\cdot ,y,z)$, we can assume,
without loss of generality, that $f(u_{1},y_{0},z)\subset f(u,y_{0},z)+C.$
We will prove that $u\in \overline{P^{V_{1}}}(x_{0},y_{0}).$ If, by
contrary, $u\notin \overline{P^{V_{1}}}(x_{0},y_{0}),$ there exists $%
z_{0}\in \overline{A^{V_{1}}}(x_{0})$ such that $f(u,y_{0},z_{0})%
\varsubsetneq C$ and, consequently, $f(u_{1},y_{0},z)\subset
f(u,y_{0},z)+C\varsubsetneq C+C\subset C,$ which contradicts $u_{1}\in 
\overline{P^{V_{1}}}(x_{0},y_{0}).$ It remains that $u\in \overline{P^{V_{1}}%
}(x_{0},y_{0}).$ Therefore, $\overline{P^{V_{1}}}(x_{0},y_{0})$ is a convex
set.

In order to prove the existence of the solutions for SVQEP, let us consider 
\ss $_{i}$ a basis of open absolutely convex symmetric neighbourhoods of
zero in $E_{i}$ for each $i\in \{1,2\}$ and let \ss =\ss $_{1}\times \text{%
\ss }_{2}.$ For each system of neighbourhoods $V=V_{1}\times V_{2}\in \text{%
\ss }_{1}\times \text{\ss }_{2},$ let's define the set valued maps $%
M^{V}:K\times D\rightarrow 2^{K\times D},$ by $%
M^{V}(x,y)=(P^{V_{1}}(x,y),F^{V_{2}}(x))=((P(x,y)+V_{1})\cap
K,(F(x)+V_{2})\cap D)$, $(x,y)\in K\times D.$ $\overline{M^{V}}$ is upper
semicontinuous with non-empty closed convex values. Therefore, according to
Ky Fan fixed point theorem \cite{fan}, there exists $(x_{V_{1}}^{\ast
},y_{V_{2}}^{\ast })\in K\times D$ such that $(x_{V_{1}}^{\ast
},y_{V_{2}}^{\ast })\in \overline{M^{V}}(x_{V_{1}}^{\ast },y_{V_{2}}^{\ast
}).$

For each $V=V_{1}\times V_{2}\in $\ss $,$ let's define $Q_{V}=\{(x,y)\in
K\times D:$ $x\in \overline{P^{V_{1}}}(x,y)\}\cap \{(x,y)\in K\times D:y\in 
\overline{F^{V_{2}}}(x,y)\}.$ $Q_{V}$ is non-empty since $(x_{V_{1}}^{\ast
},y_{V_{2}}^{\ast })\in Q_{V},$ then $Q_{V}$ is non-empty and closed. We
prove that the family $\{Q_{V}:V\in \text{\ss }\}$ has the finite
intersection property. Let $\{V^{(1)},V^{(2)},...,V^{(n)}\}$ be any finite
set of $\text{\ss\ and let }V^{(k)}=V_{1}^{(k)}\times V_{2}^{(k)},$ $%
k=1,2,...,n.$ Let $V_{1}=\underset{k=1}{\overset{n}{\cap }}V_{1}^{(k)}$ and $%
V_{2}=\underset{k=1}{\overset{n}{\cap }}V_{2}^{(k)}$. Then, $V_{1}\in $\ss $%
_{1}$ and $V_{2}\in $\ss $_{2}.$ Thus, $V=V_{1}\times V_{2}\in $\ss $%
_{1}\times $\ss $_{2}.$ Clearly $Q_{V}\subseteq \underset{k=1}{\overset{n}{%
\cap }}Q_{V^{(k)}}$ so that $\underset{k=1}{\overset{n}{\cap }}%
Q_{V^{(k)}}\neq \emptyset .$

Since $K\times D$ is compact and the family $\{Q_{V}:V\in \text{\ss }\}$ has
the finite intersection property, we have that $\cap \{Q_{V}:V\in \text{\ss }%
\}\neq \emptyset .$ Take any $(x^{\ast },y^{\ast })\in \cap \{Q_{V}:V\in $%
\ss $\},$ then for each $V\in \text{\ss },$ $(x^{\ast },y^{\ast })\in 
\overline{M^{V}}(x^{\ast },y^{\ast })$. According to Lemma 2,\emph{\ }we
have that\emph{\ } $(x^{\ast },y^{\ast })\in \overline{M}(x^{\ast },y^{\ast
}).$ Therefore, $x^{\ast }\in \overline{P}(x^{\ast },y^{\ast })$ and $%
y^{\ast }\in \overline{F}(x^{\ast }).$ We can prove by using the above
technique, that $\overline{P}$ is defined by $\overline{P}(x^{\ast },y^{\ast
})=\{u\in \overline{A}(x):f(u,y,z)\subset C,$ $\forall z\in \overline{A}%
(x)\} $ $\forall (x,y)\in K\times D.$ Consequently, there exist $x^{\ast
}\in K$ and $y^{\ast }\in \overline{F}(x^{\ast })$ such that $x^{\ast }\in 
\overline{A}(x^{\ast })$ and $f(x^{\ast },y^{\ast },x)\subset C$, $\forall
x\in \overline{A}(x^{\ast }),$ i.e. $x^{\ast }\in V_{A}(f)$ and then, $%
V_{A}(f)$ is non-empty.

\begin{remark}
Theorem 2 generalizes Theorem 3.1 in \cite{long}, since the correspondences $%
A$ and $F$ verify weaker assumptions then those ones in \cite{long}.
\end{remark}

\end{document}